\magnification=\magstephalf

\input amstex
\input epsf
\documentstyle{amsppt}

\hsize=16truecm
\vsize=21.5truecm

\def\({\left(\,}
\def\){\,\right)}
\def\[{\left[\,}
\def\]{\,\right]}

\def\<{\left\langle}
\def\>{\right\rangle}

\def\Var{\text{Var}}


\def \eqd{\buildrel d \over =}
\def \convd{\buildrel d \over \longrightarrow}

\def\a{\alpha}
\def\b{\beta}
\def\d{\delta}
\def\e{\varepsilon}
\def\g{\gamma}

\def\l{\lambda}

\def\sig{\sigma}
\def\th{\theta}

\def\z{\zeta}

\def\bC{{\bold C}}

\def\bZ{{\bold Z}}

\def\bN{{\bold N}}
\def\C{{\Cal C}}

\def\A{{\Cal A}}
\def\F{{\Cal F}}

\def\U{\,{\Cal U}}

\def\O{\Cal O}



\def\sp{\vskip1ex}
\def\z{\zeta}
\def\ps{\psi}
\def\inv{^{-1}}
\def\ra{\rightarrow}
\def\l{\ell}

\def\iy{\infty}
\def\be{\begin{equation}}
\def\ee{\end{equation}}
\def\ov{\over}
\def\al{\alpha}
\def\noi{\noindent}

\def\dl{\delta}
\def\ph{\varphi}
\def\Pr{P}
\def\rmAi{\text{Ai}}
\def\si{\sigma}

\def\eN{\eta_N}

\def\ub{\bar{u}} 
\def\lan{\left\langle} \def\ran{\right\rangle}
\def\ve{\varepsilon}



\hyphenation{per-co-la-tion}
\hyphenation{per-co-la-ting}
\hyphenation{Wads-worth}
\hyphenation{Green-berg}
\hyphenation{Has-tings}
\hyphenation{pub-li-cation}

\newdimen\howmuch

\howmuch=0pt

\def \centerbox#1{\setbox0=#1
\advance\howmuch by \the\wd0
\multiply\howmuch by -1
\advance\howmuch by \hsize
\divide\howmuch by 2
\moveright\howmuch #1}

\def\beginitems{\begingroup
    \bigbreak
    \parskip=5pt
    \advance \parindent by 2em
    \vskip-\parskip}
\def\enditems{
    \bigbreak
    \vskip-\parskip
    \endgroup}

\font\sc=cmcsc10


\document


\pageno=0

\footline={\hfil}

\null
\vskip0.5cm
 
\centerline{\bf A GROWTH MODEL IN A RANDOM ENVIRONMENT}  

\vskip1cm  

\centerline{\sc Janko Gravner}
\centerline{\rm Department of Mathematics }
\centerline{\rm University of California}
\centerline{\rm Davis, CA 95616}
\centerline{\rm email: \tt gravner\@math.ucdavis.edu}

\vskip 0.2cm

\centerline{\sc Craig A. Tracy}
\centerline{\rm Department of Mathematics}
\centerline{\rm Institute of Theoretical Dynamics}
\centerline{\rm University of California}
\centerline{\rm Davis, CA 95616}
\centerline{\rm email: \tt tracy\@itd.ucdavis.edu}

\vskip 0.2cm

\centerline{\sc Harold Widom}
\centerline{\rm Department of Mathematics}
\centerline{\rm University of California}
\centerline{\rm Santa Cruz, CA 95064}
\centerline{\rm email: \tt widom\@math.ucsc.edu}

\vskip 0.2cm

\centerline{(Version 2, July 2, 2001)}

\vskip 0.2cm

\flushpar{\sc Short Title:} Growth in random environment

\vskip 0.2cm

\flushpar{{\bf Abstract.} We consider a model of interface 
growth in two dimensions, given by a height function 
on the sites of the one--dimensional integer lattice. 
According to the discrete time update rule, the height 
above the site $x$ increases to the height above $x-1$, 
if the latter height is larger; otherwise the height above $x$ 
increases by 1 with probability $p_x$. We assume that 
$p_x$ are chosen independently at random with a common 
distribution $F$, and that the initial state is such that 
the origin is far above the other sites. We explicitly 
identify the asymptotic shape and prove that, in the 
pure regime, the fluctuations about that 
shape, normalized by the square root of time, 
are asymptotically normal. 
This contrasts with the quenched version:  
conditioned 
on the environment, and normalized by the cube 
root of time, the fluctuations almost surely approach  
a distribution known from random matrix theory. 
 
\vskip 0.2cm

\flushpar 2000 {\it Mathematics Subject Classification\/}. Primary 60K35.
Secondary 05A16, 33E17, 82B44. 

\vskip0.2cm

\flushpar {\it Keywords\/}: growth model, time constant, 
fluctuations, Fredholm determinant, Painlev{\'e} II, saddle point method.

\vskip0.2cm 

\flushpar {\bf Acknowledgments.} This work was 
partially supported by 
NSF grants DMS--9703923, DMS--9802122, and DMS--9732687,
as well as the Republic of Slovenia's Ministry of Science 
Program Group 503. We extend special thanks to 
Kurt Johansson for valuable 
insights which considerably improved 
the presentation in this paper. We also gratefully acknowledge Michael 
Casey, Bruno Nachtergaele, Timo Sepp\"al\"ainen,  
and Roger Wets for illuminating  comments.



\vfill\eject

\baselineskip=15pt

\parskip=13pt 

\pageno=1

\centerline{\bf A GROWTH MODEL IN A RANDOM ENVIRONMENT}

\vskip0.3cm  

\centerline{\sc Janko Gravner,  Craig A. Tracy, Harold Widom}

\vskip0.7cm

\subheading{1. Introduction}

Processes of random growth and deposition have a 
long history in the physics literature, 
typically as models of systems far from equilibrium (e.g., [Mea]
and the more than 1300 references listed therein). 
They made their appearance in probabilistic research 
about 35 years ago, with arguably the most basic 
growth rule, {\it first passage percolation\/} ([HW]). 
The fundamental asymptotic result is an ergodic theorem: 
scaled by time $t$, the growing 
set of sites approaches a deterministic 
limiting shape. As these early successes were 
based on nonconstructive subadditivity arguments, 
they posed two natural questions: 
(1) can the asymptotic shape be 
identified analytically and (2) how large are fluctuations about the limit?
While there has been no resolution of the first 
issue, ingenious probabilistic and geometric
arguments have yielded 
much progress on the second ([Ale]), although 
the matter is still far from settled.
It is therefore of some importance to be able to 
provide a complete answer on some other simple, but nontrivial, 
interacting growth process. 
It turns out that several two--dimensional oriented models with 
{\it last passage\/} property ([Sep1], [Sep2], [Joh1], [Joh2], [BR], 
[PS1], [GTW1]) 
are most convenient, as they can be represented, 
on the one hand, 
as particle systems related to asymmetric exclusion and,  
on the other hand, as increasing paths in random matrices 
and associated Young diagrams. This allows explicit
answers to both questions (1) and (2) above. 

In this paper we continue to study {\it Oriented 
Digital Boiling (ODB)\/} (Feb.~12, 1996, Recipe at 
[Gri], [Gra], [GTW1]), perhaps one of the 
simplest models for a coherent growing interface in the
two--dimensional lattice $\bZ^2$. 
The occupied set, which changes in 
discrete time $t=0,1,2,\dots$, is given 
by $\A_t=\{(x,y): x\in \bZ, y\le h_t(x)\}$, 
and the height function $h_t$ evolves according to the 
following rule: 
$$
h_{t+1}(x)=\max\{h_{t}(x-1), h_{t}(x)+\e_{x,t}\}. 
$$
Here $\e_{x,t}$ are independent Bernoulli random 
variables, with $P(\e_{x,t}=1)=p_x$. 
Thus the probability of a random increase depends 
on the spatial location. It remains to specify the 
initial state which will be 
$$
h_0(x)=\cases
0, &\text{if }x=0,\\
-\infty, &\text{otherwise.}
\endcases
\tag 1.1 
$$

In [GTW1] we analyzed the homogeneous case $p_x\equiv p$, 
identifying the following four asymptotic regimes: 

\beginitems

\item{--}  {\it Finite $x$ GUE Regime\/}: if $x$ is fixed and 
$t\to\infty$, then ${(h_t(x)-pt)}/{\sqrt{p(1-p)t}}\convd M_x,$ 
a Brownian functional  
whose law can be computed explicitly 
as the largest eigenvalue of $(x+1) \times (x+1)$ hermitian
matrix from the Gaussian Unitary Ensemble (GUE).

\item{--} {\it GUE Universal Regime\/}: if $x$ is a positive multiple of $t$,
and $x/t<1-p$, then there exist constants $c_1$ and $c_2$ so that 
${(h_t(x)-c_1t)}/{(c_2 t^{1/3})}$ converges weakly 
to a distribution $F_2$ known from random matrix theory ([TW1]). 
 
\item{--}  {\it Critical Regime\/}: if $x=(1-p)t+o(\sqrt t)$, then 
$
P(h_t(x)-(t-x)\le -k)
$
converges to a $k\times k$ determinant.

\item{--}  {\it Deterministic Regime\/}: 
if $x$ is a positive multiple of $t$, and $x/t>1-p$, 
$P(h_t(x)=t-x)\to 1$ exponentially fast.

\enditems

The focus of this paper is ODB
in a random environment, in which $p_x$ are 
initially chosen at random, with common
distribution given by $P(p_x\le s)=F(s)$. We will also 
assume that $p_x$ are  independent, although 
in several instances this assumption can be  
considerably weakened. 
In statistical physics, processes in a random environment 
are often called {\it disordered systems \/}, or, especially
in the  Ising--type models, {\it spin glasses\/}. 
In this context, the random environment (choice of $p_x$) 
is referred to as {\it quenched\/} randomness, as opposed to
the dynamic ({\it thermal\/}) fluctuations 
induced by the coin flips $\e_{x,t}$.  
In general, 
rigorous research in this area has been a notoriously difficult 
enterprise; for some recent breakthroughs (as well as reviews of 
the literature) we refer the reader to  [SK], [NS], [NV] and [Tal]. 

We now state our main results. Throughout, 
we will denote by $\<\,\cdot\,\>$ integration with 
respect to $dF$ and $p$ a generic random variable 
with distribution $F$.  

Construct a random $m\times n$ matrix $A=A(F)$, with 
independent Bernoulli entries $\e_{i,j}$ and such 
that $P(\e_{i,j}=1)=p_j$, where, again, $p_j\eqd p$ are i.i.d.
Label columns as usual, 
but rows started at the bottom. We call a sequence 
of 1's in $A$ whose positions have column index nondecreasing and row index 
strictly increasing an {\it increasing path\/} in $A$. 
Let $H=H(m,n)$ be the length of the longest  
increasing path. (Sometimes, to emphasize 
dependence on $F$, we write $H=H(F)=H(m,n,F)$.) The following lemma is then 
easy to prove ([GTW1]). 

\proclaim{Lemma1.1} Under a simple coupling, $h_t(x)=H(t-x,x+1)$. 
\endproclaim

We will therefore concentrate our attention on the random 
matrix $A$ from now on, switching to the height function 
only occasionally to interpret the results. We also note 
that Lemma 1.1 demonstrates that ODB is 
equivalent to the Sepp\"al\"ainen--Johansson model 
([Sep2], [Joh2]).  

Our first theorem identifies the time constant. In the 
sequel, we will present two completely different 
methods for proving these limits, a variational approach 
and a determinantal approach. The first 
method (which is similar to the one in [DZ]) 
is based on the crucial symmetry 
property of $H$ (Lemma 2.2) and provides some 
information on the longest increasing path itself,  
while the second one is deeper and more precise 
and thus able also to determine fluctuations. 
The paper [SK] studies a
related model, presents 
yet another technique, based on an exclusion process 
representation, and observes similar 
phase transitions. Throughout this paper, we 
let 
$$b=b(F)=\min\{s:F(s)=1\}$$
be the  
right edge of the support of $dF$ and assume 
that $n=\a m$ for some $0<\a<\infty$. (Actually, 
$n=\lfloor \a m\rfloor$, but we  
drop the integer part as it obvious where 
it should be used and to avoid complicating expressions.) We also 
define the following critical values 
$$
\aligned
&\a_c=\< \frac{p}{1-p}\>^{-1},\\
&\a_c'=\<\frac{p(1-p)}{(b-p)^2}\>^{-1}
\endaligned 
\tag 1.2
$$
and define 
$c=c(\a,F)$ to be the time constant 
$$
c=c(\a,F)=\lim_{m\to\infty}\frac{H}m.
\tag 1.3
$$
Note that $c$ determines the limiting shape 
of $\A_t$, namely $\lim \A_t/t$, as $t\to\infty$
for the corner initialization given by (1.1). 
By virtue of the Wulff transform, it then also 
gives the speeds of {\it some\/} half--planes, i.e., 
$\lim \A_t/t$ when $\A_0$ comprises points 
below a fixed line. See [SK] for much more on this 
issue. 

\proclaim{Theorem 1} The limit in (1.3) exists 
almost surely. If $b=1$, then $c(\a,F)=1$ for all $\a$, 
while if 
$b<1$, then 
$$
c(\a,F)=
\cases
b+\a (1-b)\<p/(b-p)\>,&\text{ if }\a\le \a_c', \\
a+\a (1-a)\<p/(a-p)\>,&\text{ if }\a_c'\le \a\le \a_c,\\
1,&\text{ if }\a_c\le \a. 
\endcases
$$
Here $a=a(\a,F)\in [b,1]$  is the unique solution to
$$\a\<\frac {p(1-p)}{(a-p)^2}\>=1.$$
\endproclaim

Note that that $\<(b-p)^{-2}\>=\infty$ 
iff $\a_c'=0$ iff there is only one 
critical value. 

Next we turn our attention to fluctuations. In this paper we  
present complete results for the {\it pure\/} regime $\a_c'<\a< \a_c$
and for the (easy) {\it deterministic\/} regime $\a_c<\a$. 
The {\it composite\/} regime $\a< \a_c'$ is addressed 
in [GTW2], while both {\it critical\/} cases when 
$\a$ equals either critical value currently remain unresolved. 
To explain the results, and connect with the spinglass 
terminology we have just used, we turn to a simulation. 
For an example, we use 
$F(s)=1-(1-2s)^{3}$  
so that $b=1/2$, $\a_c\approx 6.3$ and 
$\a_c'\approx 0.5$ and run the simulation until time $t=40,000$
(with a single realization of the environment and the 
coin flips).  When 
$x$ is close to the origin, it is clear from the 
picture that the interface mostly consists of sheer walls
followed by flat pieces. The walls correspond to the rare 
sites with update probability $p_x$ close to $1/2$. Those 
are much faster than the other sites so they pull 
ahead of their left neighbors, creating walls, and 
dominate their right neighbors by ``feeding'' them
at nearly largest possible rate. In fact, this state 
of affairs persists up to about $x=t/3$ although close to $x=t/3$ 
these effects are less pronounced. In the pure regime, when 
$x/t$ ranges approximately from  $0.333$ to approximately $0.863$, 
the fluctuations are much more regular, and in fact, as 
we will demonstrate, asymptotically normal. For larger $x/t$ the shape 
has slope $-1$ and no fluctuations.

\smallskip
\epsfysize=0.19\vsize
\centerline{\epsffile{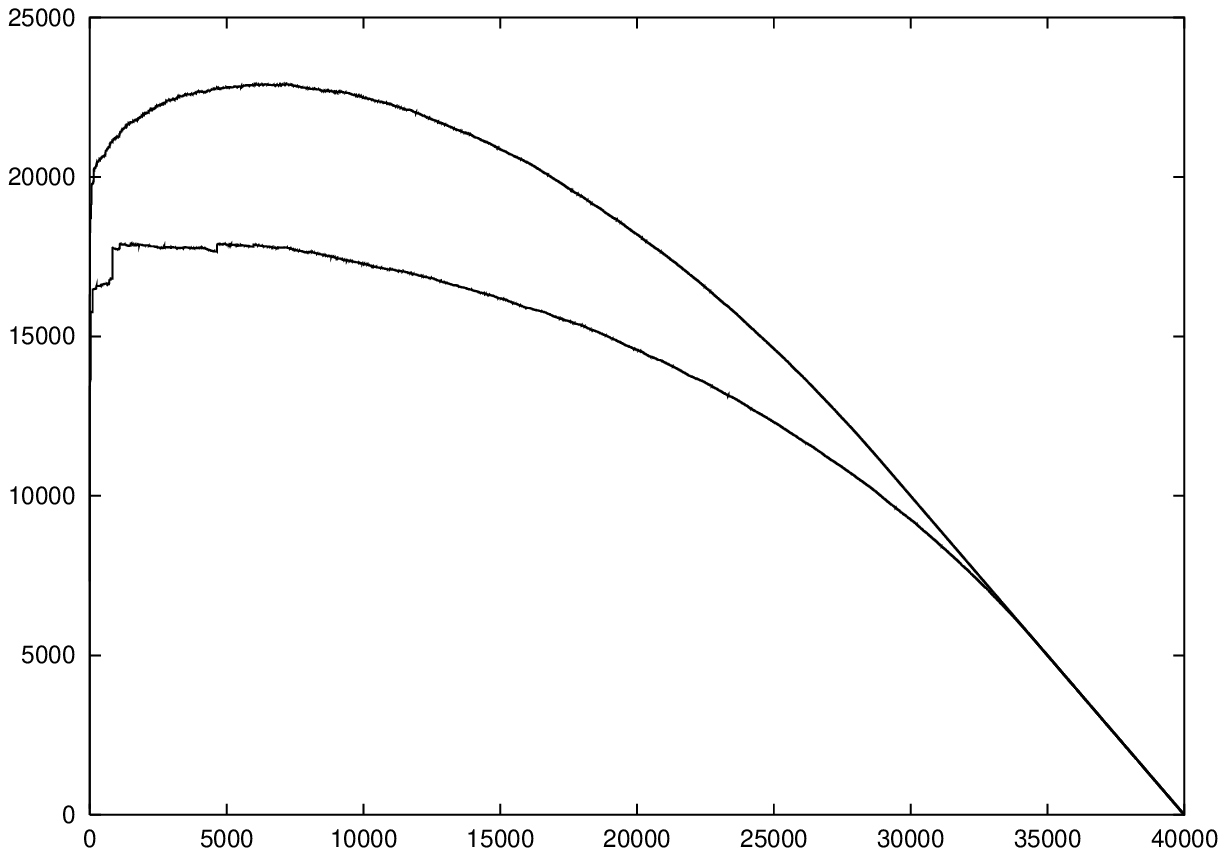}}
\smallskip
\epsfysize=0.19\vsize
\centerline{\epsffile{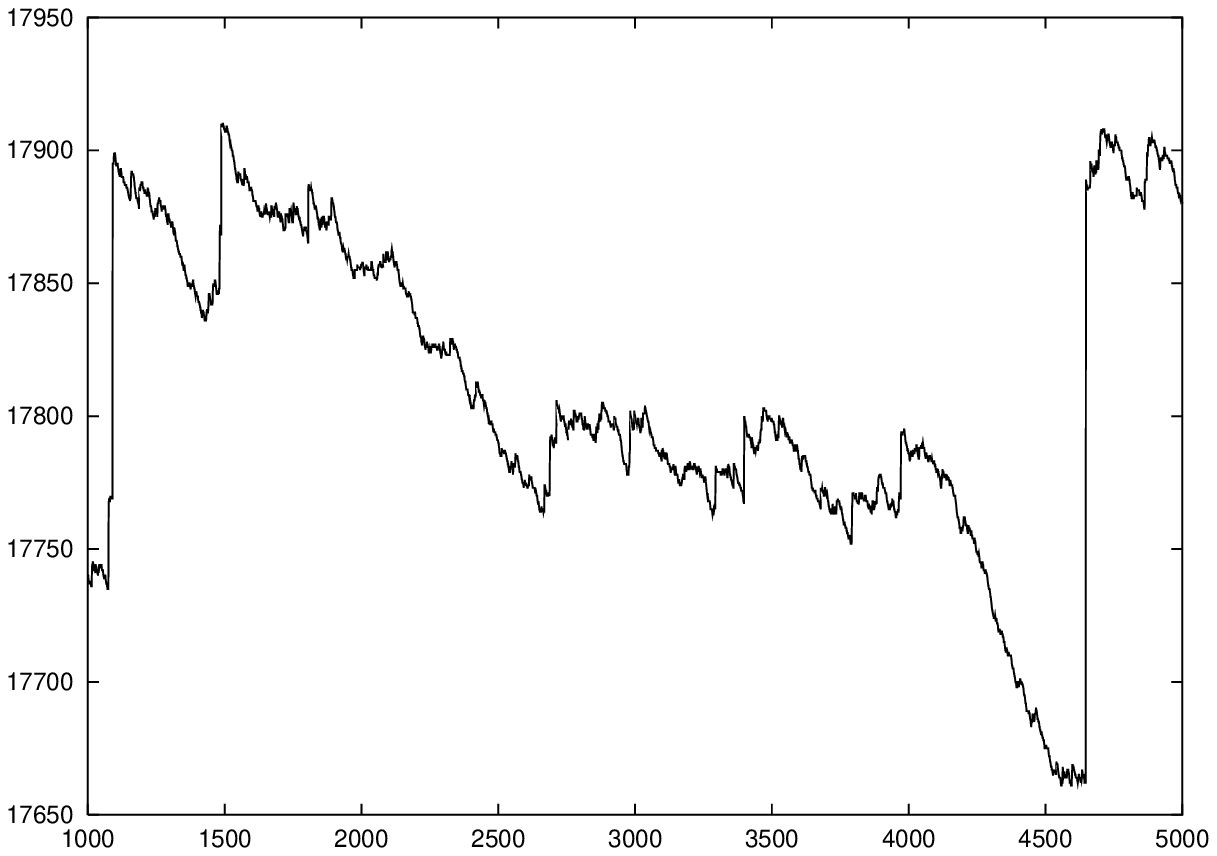}
\epsfysize=0.19\vsize
\epsffile{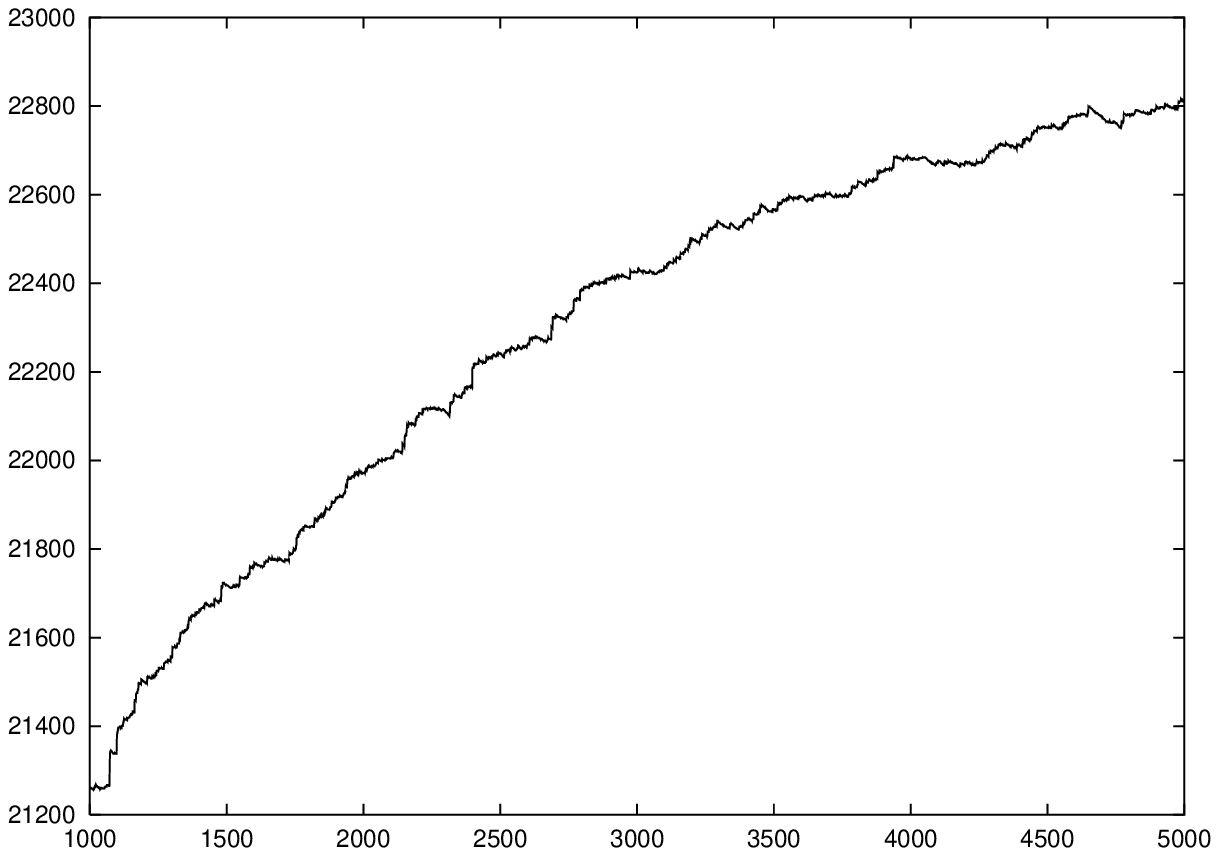}}
\smallskip
\centerline {Figure 1. Two ODB simulations, as explained in the text.}
\smallskip

For comparison, consider the case when $p$ is uniform on 
$[0,1/2]$, the case that has 
$\a_c=1/(\ln 4-1)\approx 2.59$ 
and $\a_c'=0$. The fluctuations are normal up to $x/t\approx 0.72$. 
Figure 1 depicts the results of simulations, first 
complete boundaries of two occupied sets (the 
top curve is the uniform case), then two details (the right curve is 
the uniform case)
for $x\in [1000, 5000]$.

\proclaim{Theorem 2} Assume that $b<1$ and $\a_c'<\a<\a_c$. 
Let $a$ be as in Theorem 1 and 
$$
\tau^2=\Var\(\frac{(1-a)p}{a-p}\). 
$$
Then, as $m\to\infty$, 
$$
\frac{H-cm}{\tau\sqrt\a \cdot m^{1/2}}\convd N(0,1).
$$
\endproclaim

Assume that $p$ is uniform $[0,1/2]$ to 
illustrate Theorems 1 and 2. Together 
they imply that there exist $c_1$ and $c_2$ so 
that $(h_t(x)-c_1t)/(c_2 t)^{1/2}\convd N(0,1)$,
where  
$c_1$ determines the limiting shape and $c_2$ is the variance. 
These two quantities are presented in Figure 2, $c_1$ is the 
top and $c_2$ is the bottom curve.  For comparison, the 
shape of homogeneous ODB with $p_x\equiv \<p\>=1/4$
is also drawn (middle curve). 
Note that $c_1$ and $c_2$ approach $1/2$ and $1/4$, respectively, 
as $\a\to 0$, indicating that for small $x/t$ the interface 
growth is governed by the largest update probability, 
which is close to 1/2.  
Finally, we do the same computation for the other example 
in Figure 1. The variance is now drawn only on $[\a_c', 1]$.

\smallskip

\epsfysize=0.19\vsize
\centerline{\epsffile{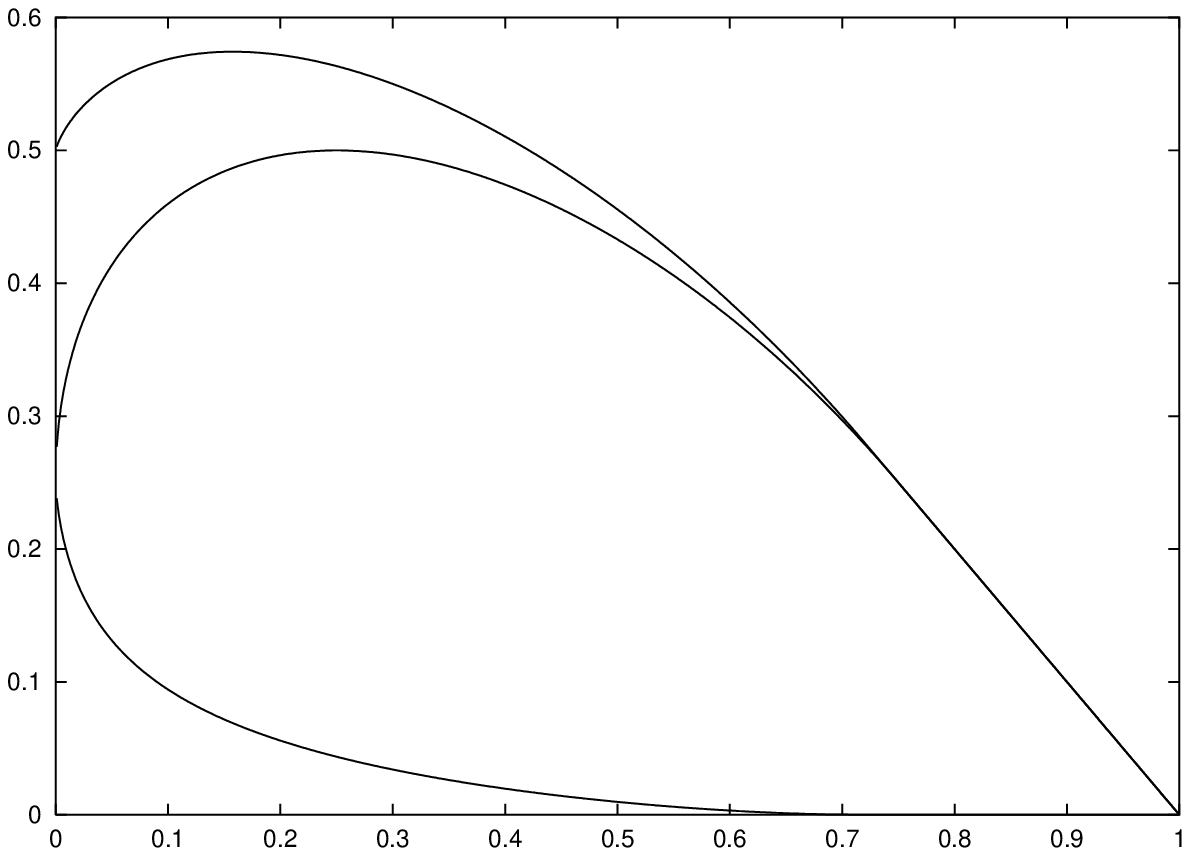}
\epsfysize=0.19\vsize
\epsffile{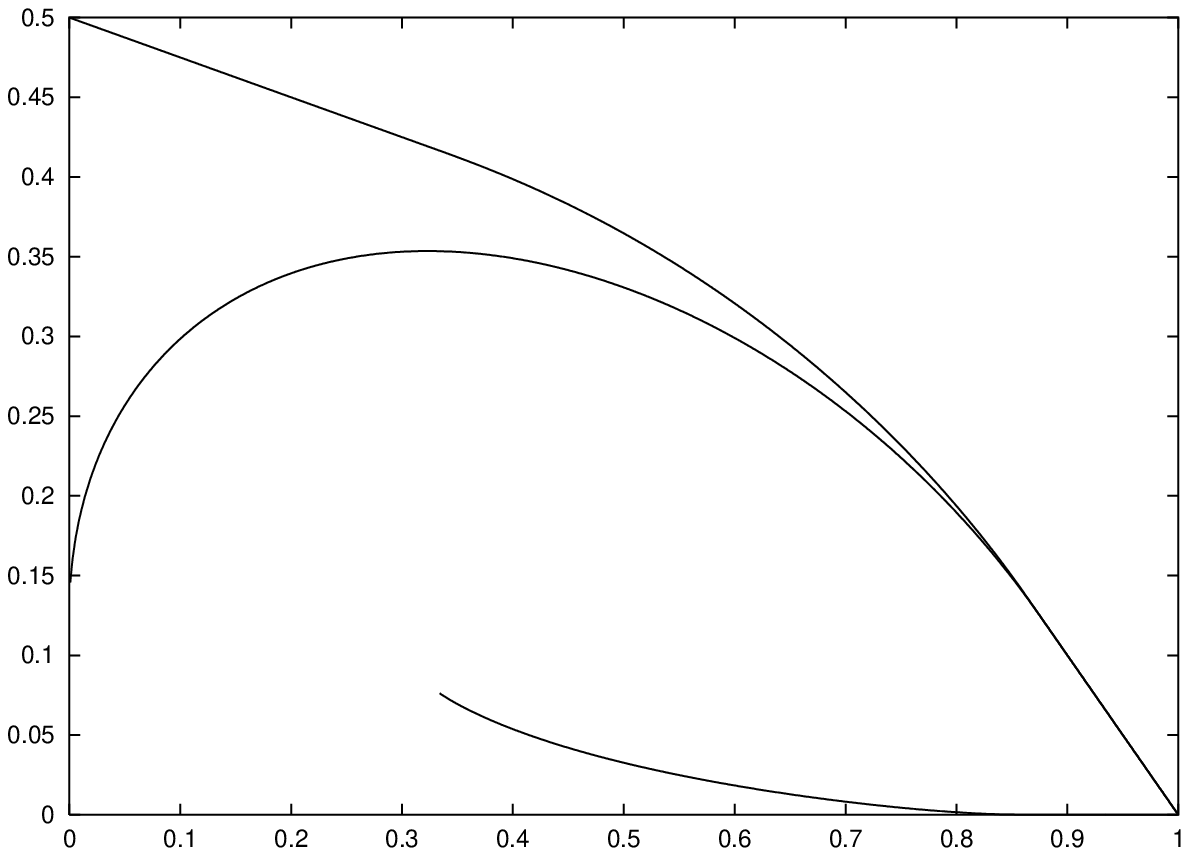}}


\smallskip
\centerline {Figure 2. $c_1$ (top), $c_2$ (bottom)
and the shape for $p_x\equiv \<p\>$ (middle) vs. $x/t$.}
\centerline {The two distributions are uniform [0,1/2] (left)
and $F(s)=1-(1-2s)^{3}$.}

\smallskip

\comment

\smallskip
\epsfysize=0.25\vsize
\centerline{ex3.ps}
\smallskip
\centerline {Figure 3. Same quantities as in Figure 2 for the case 
$F(s)=1-(1-2s)^{3}$.}
\smallskip

\endcomment

We note that both a.s.\ convergence to the limiting
shape (which is equivalent to a.s.\ convergence in (1.3)) 
and its  convexity follow from subadditivity, which in turn 
is a consequence of the fact that this is
an oriented model in which influences only travel in 
one direction. To be more precise, fix integer sites 
$(x_1,y_1), (x_2,y_2)\in \bZ_+\times \bZ_+$  
and define times $T_{(x_1,y_1), (x_2,y_2)}$ as follows. 
First wait until time
$T_{(0,0),(x_1,y_1)}$ when   
the dynamics reaches $(x_1, y_1)$. 
Then restart the dynamics from the initial state 
$$
h_0(x)=\cases
y_1, &\text{if }x=x_1,\\
-\infty, &\text{otherwise.}
\endcases
$$
and let $T_{(x_1,y_1), (x_2,y_2)}$ be the time at which the occupied 
set reaches $(x_2,y_2)$. This random variable is 
independent of $p_x$ for $x\le x_1-1$ and 
$T_{(0,0),(x_2,y_2)}\le T_{(0,0),(x_1,y_1)}+ T_{(x_1,y_1), (x_2,y_2)}$.
Therefore, the subadditive ergodic theorem can be applied as
in the first chapter of [Dur]. 

The main step in the proof of Theorem 2 establishes a 
limit 
law for fluctuations conditioned on the state of the
environment. In many ways, such a result is more 
pertinent to understanding physical processes 
modeled by simple growth models such as ODB.

\proclaim{Theorem 3} Assume that $b<1$ and 
$\a_c'<\a<\a_c$. Then there exists a sequence of
random variables
$G_n\in \sigma\{p_1,\dots, p_n\}$ and a constant $g_0\ne 0$ 
(both depending on $\a$)
such that, as $m\to \infty$,   
$$
P\(\frac{H-G_n}{g_0^{-1} m^{1/3}}\le s\,\mid\, p_1,\dots,p_n\)\to F_2(s), 
$$
almost surely, for any fixed $s$. 
\endproclaim

The random variables $G_n=c_nm$ 
are given in terms of the solution 
of an algebraic equation in which $p_1,\dots, p_n$ 
appear as parameters (see (3.4) and (3.5)), while the deterministic constant
$g_0$ is specified before the statement of Lemma 3.5. 
The limiting distribution function $F_2$ first arose
in connection with eigenvalues 
of random matrices ([TW1], see [TW2] for a review). Since then
it has been observed in many other contexts, including 
growth processes ([Joh1], [Joh2], [BR], [GTW1], [PS1], [PS2]). Most suitable 
for computations is the identity
$$
F_2(s)=\exp\left(-\int_s^\infty
(x-s) q(x)^2\,dx\right),  
$$
where $q$ is the unique solution of the Painlev\'e II equation
$$q''=sq+2q^3,$$
which is asymptotic to the Airy 
function, $q(s)\sim\text{Ai}(s)$ as $s\ra\infty$. 
When proving limit laws, it is more 
useful that  
$F_2$ can be represented as a Fredholm determinant
(see e.g. [GTW1] and Section 3 below). 

In Theorem 3 the environment is assumed as given, 
$H$ is approximated by the quenched shape $G_n$, 
with the fluctuations about this shape of
the order $m^{1/3}$ and given by the $F_2$ distribution. 
As we prove in 
Section 3, $(\a m)^{-1/2}(G_n-c m)$ converges to the standard normal, 
making it clear why Theorem 2 holds: the environmental 
noise eventually drowns out the more interesting quenched 
fluctuations of Theorem 3. An illustration is provided 
in Figure 3, in which $p$ is again uniform on $[0,1/2]$ 
and $h_t$ (solid curve), deterministic approximation 
based on Theorem 2 (dotted curve), and the much 
better random approximation 
based on Theorem 3 (dashed curve) are all depicted at times
$t=100, 200, \dots, 1000$.

\smallskip
\epsfysize=0.35\vsize
\centerline{\epsffile{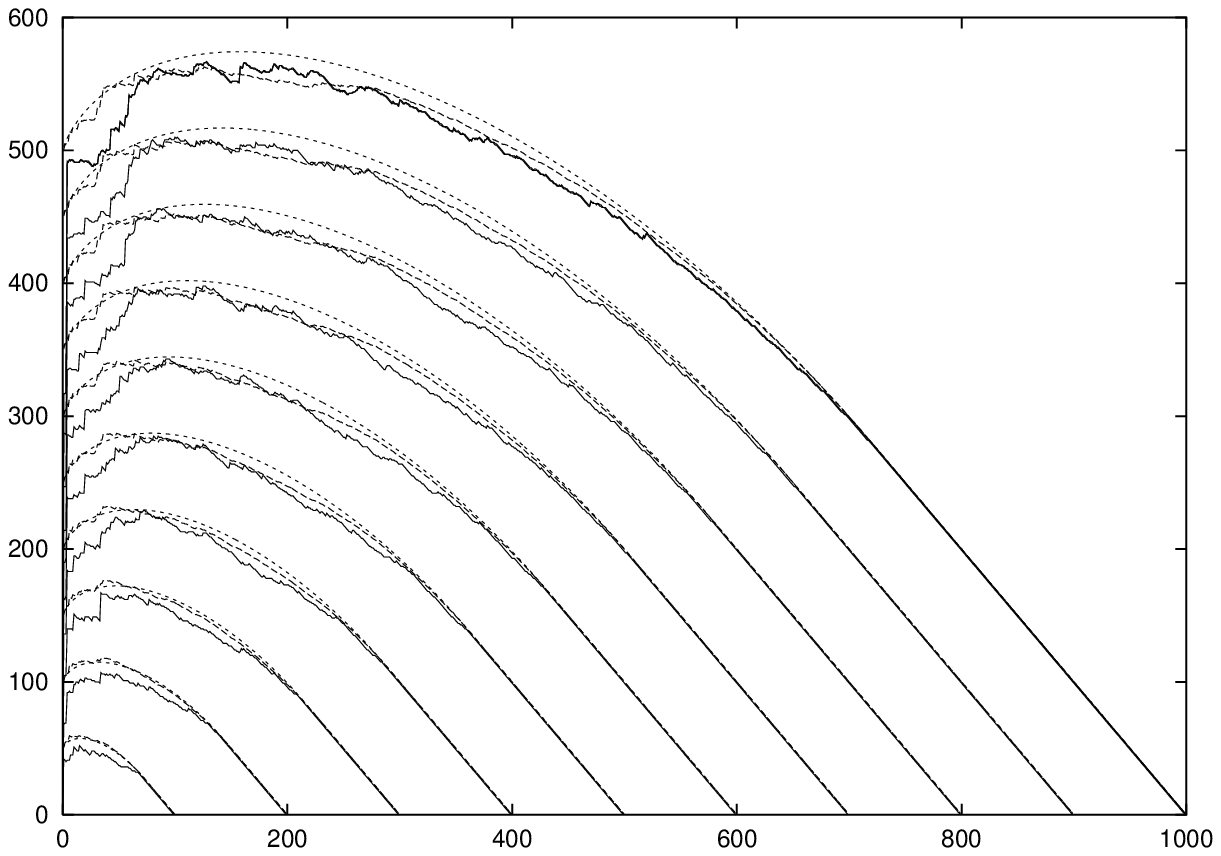}}
\smallskip
\centerline {Figure 3. Approximations to $h_t$ (solid curve) based on}
\centerline{Theorems 2 (dotted curve) and 3 (dashed curve).}

In conclusion, we note that the 
connection between random matrix theory and 
random combinatorial objects, 
which has become the key to rigorous 
understanding of random interface fluctuations, 
made its initial appearance in [BDJ], while
an inhomogeneous model
of ODB type was first studied in [ITW1]. This last paper, 
together with its companion [ITW2],  
extends the study of random words from the 
homogeneous case in [TW3] in a somewhat analogous 
way as the present paper builds on the work in [GTW1]. 
In particular, connections with operator determinants 
(from the beginning of Section 3) are very similar (see 
also [Rai] which features a general inhomogeneous setup). However,
randomness of the environment, 
which seems to be a new feature in rigorous 
analysis of explicitly solvable 
models, 
then forces our techniques to take a novel turn.

\subheading{2. A variational characterization of the time constant}

We start by a remark on constructing the random matrix $A$. 
The most convenient design uses as the 
probability space $(\Omega, P)$ a countably infinite product of  
unit intervals $[0,1]$ with Lebesgue measure. A copy 
of the unit interval (and thus a factor in the 
product) is associated with each point 
in  
$\bN\times \bN$  
and in addition, with each 
positive integer in $\bN$. 
(The former factors correspond to matrix entries, and 
the latter to its columns.) If $\omega=(m_{ij}, c_j)\in \Omega$ 
is a generic realization, we define the following 
random variables: 
$
p_j=F^{-1}(c_j)
$
(where $F^{-1}(x)=\sup\{y:F(y)<x\}$ as usual) and 
$\e_{ij}=1_{\{m_{ij}<p_j\}}$. By restricting to the  
$m\times n$ rectangle at the lower right corner of 
$\bN\times \bN$, this  
constructs the random matrices $A$ for 
all $m$ and $n$ simultaneously.  
The following useful lemma
also follows immediately. 

\proclaim{Lemma 2.1} If $F_1\le F_2$ are two distribution 
functions, the two corresponding random matrices $A(F_1)$ and 
$A(F_2)$  
can be coupled so that $H(F_2)\le H(F_1)$. 
\endproclaim 

Next we state the crucial property for the variational  
approach to work: conditioned on the environment, 
$H$ is a symmetric function of flip probabilities. 

\proclaim{Lemma 2.2} A regular conditional distribution 
$$
P(H\le h\,|\, p_1,\dots, p_n)
$$
is a symmetric function of $p_1,\dots p_n$. 
\endproclaim

\demo{Proof}
See section 2.2 of [GTW1]. $\square$
\enddemo

Somewhat loosely, we denote by $H_n$ the random variable 
$H$ obtained by fixing $p_1,\dots, p_n$. 
In fact this is nothing more 
that a shorthand notation, e.g., 
$E(\varphi(H_n))=E(\varphi(H)\mid p_1,\dots, p_n)$
for any bounded measurable function $\varphi$.

The time constant $c(\a,x)=c(\a,\d_x)$ for the case 
$p_j\equiv x$ is given in [GTW1]. The next lemma summarizes 
the relevant conclusions. 

\proclaim{Lemma 2.3} Assume that $dF=\d_x$. Then  
$$
c=c(\a,x)=\cases
2\sqrt\a\sqrt{x(1-x)}+(1-\a)x,\quad&(1-x)/x>\a,\\
1,\quad&(1-x)/x\le\a.
\endcases
$$
Moreover, for every $\e>0$ there exists a constant $\g=\g(\e)>0$ 
so that  
$$
P(|H/m-c|>\e)<e^{-\g m} 
\tag 2.1 
$$
for $m\ge m_0(\e,\a,x)$. 
\endproclaim

\demo{Proof} The formula for $c$ follows from (3.1) 
in [GTW1], while the large deviation estimate can 
be proved by the method of bounded differences 
as in Lemma 5.4 of [Gra]. $\square$
\enddemo

It turns out the following function is more 
convenient than $c$. 
$$
\z(y,x)= y \cdot c(1/y,x)=\cases
2\sqrt y\sqrt{x(1-x)}+(y-1)x,\quad&x/(1-x)< y,\\
y,\quad&x/(1-x)\ge y.  
\endcases
$$
Note that the partial derivative  
$$
\z_y(y,x)=
\cases
y^{-1/2}\sqrt{x(1-x)}+ x,\quad&x/(1-x)< y ,\\
1, \quad&x/(1-x)\ge y .
\endcases
$$
is decreasing in $y$ (obviously) and increasing in $x$
(easily checked). In particular, $\z(\cdot,x)$ is a convex
function.  

We now  derive a variational 
problem for $c$, initially without paying  attention 
to rigor. Start by a nice distribution function $F$ and  
approximate it by the discrete distribution 
function given by 
$$
P\(p_j=\frac ik\)=
\Delta F_k(i)=F\(\frac{i}k\)-F\(\frac{i-1}k\),
 i=1, \dots, k. 
$$
Let $\psi:[0,\a]\to[0,1]$, $\psi(0)=0$, 
$\psi(\a)=1$ be a nondecreasing 
function, with $\Delta\psi_k(i)= \psi(\a F(i/k))-\psi(\a F((i-1)/k))$. 
Define the functionals: 
$$
\F(\psi)=\int_0^1 \z(\psi'(\a F(x)), x)\cdot \a\, dF(x) 
$$
and 
$$
\F_k(\psi)=
\sum_{i=1}^k 
\z\(\frac{\Delta\psi_k(i)}{\a\Delta F_k(i)}, \frac ik\)
\cdot  {\a\Delta F_k(i)}.
$$

Generate the 
$p_j$'s and denote by 
$N_i$ 
the number of $p_j$ equal to $i/k$.
By Lemma 2.2,  we can assume flip probability $1/k$ 
in the first $N_1$ columns, 
$2/k$ in the next $N_2$ columns, etc.
Moreover, the strong law suggests that
the identity $N_i=\Delta F_k(i)n$ nearly holds. 
As we know the asymptotics 
for the
longest increasing paths in the slivers of widths $N_i$ 
in which the probabilities are constant, the longest 
increasing path in $A$ is determined by the 
most advantageous choice of transition 
points between the slivers. These 
transition points are specified by 
a function $\psi$ as described above. If we  
approximate the differences with 
derivatives, we obtain
$$
\align 
&c(\a, F)=\lim_{k\to\infty}c(\a,F_k)\\
&=\lim_{k\to\infty}\max_{\psi}\sum_{i=1}^k 
c\(\frac{\a\Delta F_k(i)}{\Delta\psi_k(i)}, \frac ik\)
\cdot \Delta\psi_k(i)\\
&=\lim_{k\to\infty}\max_{\psi} \F_k(\psi)\\
& =\lim_{k\to\infty} \max_{\psi}\sum_{i=1}^k 
\z\(\psi'(\a F(i/k)), \frac ik\)\cdot 
\a F'(i/k)\cdot\frac 1k\\
&=\max_{\psi} \F(\psi). 
\endalign
$$

At this point, we remark that a connection between 
longest increasing paths and variational problems 
has appeared before in the literature. The result closest 
to ours is by  
Deuschel and Zeitouni ([DZ]), who used a variational approach 
to study a variant of Ulam's problem.  
In their case, a number of points in the unit square is chosen 
independently according to some distribution with 
a density, then a longest sequence, increasing in both coordinates, 
is extracted from this sample. The Deutchel--Zeitouni 
functional is different from ours as the length of the 
longest increasing path has a non--trivial dependence on $\a$
(that is, through $c$) in our case. 

The (integrated) Euler functional for the variational problem is 
$$
\z_y(\psi'(x),F^{-1}(x/\a))=a, 
$$
or, writing $g(x)=\psi'(\a F(x))$,  
$$
\z_y(g(x), x)=a. 
\tag 2.2
$$
Since $\z_y\le 1$ and equal to $1$ if and only if
$x/(1-x)\ge y$,
the integration constant $a\in [0,1]$. If $a=1$, then 
$g(x)\le x/(1-x)$, $\z(g(x), x)=g(x)$ and 
$$
c(\a,F)=\int_0^1\psi'(\a F(x))\a\,dF(x)=1.
$$
Assume now that $b<1$. In this case, it is necessary to specify 
$g$ only on $[0,b)$. But (2.2) gives 
$$
g(x)=\frac{x(1-x)}{(a-x)^2}.
\tag 2.3 
$$ 
The constant $a$ is given by the boundary conditions. 
Assuming that (2.3) holds on $[0,b]$,  
$$
1=\a\int_0^b g(x)\, dF(x)=\a\int_0^b \frac{x(1-x)}{(a-x)^2}\, dF(x).
\tag 2.4
$$
The smallest the last integral can be is when $a=1$, which 
yields the condition 
$$
1>\a\int_0^b \frac{x}{1-x}\, dF(x)=\frac{\a}{\a_c}.   
$$
On the other hand, the largest that the integral in (2.4)  
can be is when $a=b$. Therefore,  if $\a\in(\a_c',\a_c)$, 
we have found the minimizer and 
$$
c(\a,F)=\int_0^b\z (g(x), x)\a \, dF(x)=\a\<\frac {-p^2-a^2 p+2ap}
{(a-p)^2}\>,
$$
which reduces, upon using the defining equation for $a$,
to the formula in Theorem 1. 

If $\a<\a_c'$, the minimizer $\psi$ 
has to make a jump of size $1-\a/\a_c'$ at $\a$. The natural 
interpretation for this is that the minimizer given by (2.3) is 
used in the lower left part of $A$ 
with dimensions $(\a/\a_c')m\times (n-1)$.  
To the resulting increasing path in this 
submatrix one needs to add the number of 1's in 
the upper segment of length $(1-\a/\a_c')m$ in the last column, 
in which nearly the largest probability $b$ is used. Therefore, 
$$
c(\a, F)=c(\a_c', F)\cdot\frac{\a}{\a_c'}+b\(1-\frac{\a}{\a_c'}\), 
$$
which again reduces to the appropriate formula in Theorem 1. 

We now proceed to give a proof Theorem 1, the heart of which
is a somewhat involved multistage approximation scheme.

\demo{Proof of Theorem 1 when $b=1$}
This follows simply by observing that, for any
$\e>0$, $\max_j p_j\to b$  a.s.\ as $m\to\infty$. 
Since a trivial lower bound is obtained 
by using only the column with 
the largest $p_j$, one concludes
that $\liminf H/m\ge b$ a.s. 
$\square$
\enddemo 

\demo{Proof of Theorem 1 when $\a\in (\a_c',\a_c)$}

We begin by the following lemma. 

\proclaim{Lemma 2.4} Assume that a sequence of distribution 
functions $F_N$ converges to $F$ in the usual sense 
(i.e., the induced measures
converge weakly). Assume also that $b(F_N)\to b(F)$ and 
that $\a_N\to \a$. Then $c(\a_N,F_N)\to c(\a,F)$ 
(as given in Theorem 1). 
\endproclaim

\demo{Proof} If $a'>b(F)$ and 
$$\a\int x(1-x)(a'-x)^{-2}\, dF(x)>1,$$ 
then for a large $N$, $a'>b(F_N)$ and,  
since the integrand is bounded, 
$$\a_N\int x(1-x)(a'-x)^{-2}\, dF_N(x)>1.$$
Hence $a_N=a(\a_N, F_N)>a'$. If $a_N\to a_0$, then 
$x(1-x)(a_N-x)^{-2}$ converges to $x(1-x)(a_0-x)^{-2}$ 
uniformly for $x\in [0,a']$ and so 
$$
1=\a_N\int x(1-x)(a_N-x)^{-2}\, dF_N(x)\to 
\a\int x(1-x)(a_0-x)^{-2}\, dF(x). 
$$
Therefore $a_0=a(\a,F)$ and consequently $a_N\to a(\a,F)$. 
As $x(a_N-x)^{-1}$ also converges uniformly on $[0,a']$, 
$$
\align 
c(\a_N, F_N)&=a_N+\a_N(1-a_N)\int x(a_N-x)^{-1}\, dF_N\\
&\to 
a+\a(1-a)\int x(a-x)^{-1}\, dF=c(\a,F).
\endalign 
$$
$\square$
\enddemo

First we assume that $F$ is {\it nice\/}, that is, a one-to-one 
function on $[\b,b]\subset (0,1)$, 
with $F(\b)=0$, $F(b)=1$, and continuously 
differentiable on $(0,1)$. We also assume that 
$\Psi$ is the class of non--decreasing convex functions 
$\psi\in \C^2[0, \alpha]$, with $\psi(0)=0$, $\psi(\a)=1$, 
$\psi'(0)\ge\b/2$. This last assumption is necessary because
$\z(y,x)$ is not Lipshitz near $y=0$.   
 
\proclaim{Lemma 2.5} Assume that $\a\in (\a_c',\a_c)$. 
Among all $\psi\in \Psi$, the 
functional $\F(\psi)$ is uniquely maximimized by 
$$
\psi(x)=\int_0^x g(F^{-1}(u/\a))^2 \, du, 
$$
where $g$ is given by (2.3).   

\endproclaim

\demo{Proof} This follows from standard calculus of 
variations. Both $\psi'(0)\ge\b/2$ 
and convexity of $\psi$ are easily checked. 
$\square$
\enddemo

We now justify the approximation 
steps in the heuristic argument, using the 
same notation. 
First, if $\e>0$ is fixed, then with 
probability exponentially (in $n$) close to 1, 
$$
(1-\e)\Delta F_k(i)n\le N_i\le (1+\e)\Delta F_k(i)n
$$
for every $i=1,\dots, k$. By obvious monotonicity, 
the longest increasing path in $A$ is then bounded 
above by the longest increasing path in $A'$
in which all $N_i=(1+\e)\Delta F_k(i)n$, and therefore we 
can get an upper bound by  increasing 
$\a$ to $\a(1+2\e)$ and assuming $N_i=\Delta F_k(i)n$. 
A lower bound is obtained similarly. As our final characterization 
of $c$  
is continuous with respect to $\a$ (Lemma 2.4), 
we can, and will, assume that $N_i=\Delta F_k(i)n$
from now on. 

The above paragraph eliminates randomness of $p_j$'s;  
we now proceed to replace the coin flips with 
deterministic quantities. Again, fix an $\e>0$ and let $M=\e m$.  
For $j_1\le j_2$ and $i=1,\dots, n$, consider 
the longest increasing paths $\pi_{j_1,j_2,i}$ 
between $(F_k(i-1)n,j_1)$ 
(noninclusive) 
and $(F_k(i)n, j_2)$ (inclusive). 
Then, 
with probability exponentially close to 1, 
the length of any  $\pi_{j_1,j_2,i}$ is at most 
$$
(1+\e) c\(\frac {\Delta F_k(i)n}{M\lceil (j_2-j_1)/M\rceil},\frac ik\)
\cdot M\lceil (j_2-j_1)/M\rceil
=(1+\e)\z\(\frac{M\lceil (j_2-j_1)/M\rceil}{\Delta F_k(i)n}\)
\cdot\Delta F_k(i)n. 
$$
(This uses Lemma 2.3 when $j_2-j_1$ is divisible by $M$ and fills
the rest by monotonicity. Note that Lemma 2.3 is therefore 
only applied finitely many times for fixed $\e$ and $k$.) 
The lower bound is obtained by rounding down instead of up. 
It follows that the length of any $\pi_{j_1,j_2,i}$
is bounded above (resp. below) by 
$$
\z\(\frac{(j_2-j_1)}{\Delta F_k(i)n}\)\cdot\Delta F_k(i)n. 
\tag 2.5
$$
computed on the matrix of size $(m+M)\times n$ (resp. $(m-M)\times n$). 
Once again we can use continuity to assume that the length 
of any $\pi_{j_1,j_2,i}$ is given by (2.5). 

It remains to show that the 
discrete deterministic optimization problem $\max_\psi \F_k(\psi)$ 
is for large $k$ close to its continuous counterpart $\max_\psi \F(\psi)$.  
To this end, we first prove that we can indeed restrict the 
set of function $\psi$ to those in $\Psi$, i.e., those
that are convex and have a large enough derivative.  
Let $\Delta x_1=\alpha\Delta F_k(i)$, 
$\Delta x_2=\alpha\Delta F_k(i+1)$,  $\Delta y_1=\Delta \psi_k(i)$,
$\Delta y_2=\Delta \psi_k(i+1)$, $\Delta y=\Delta y_1+\Delta y_2$, 
$p_1=i/n$, $p_2=(i+1)/n$. Then 
$$
\z\(\frac {\Delta y_1}{\Delta x_1}, p_1\)\Delta x_1 
+\z\(\frac {\Delta y-\Delta y_1}{\Delta x_2}, p_2\)\Delta x_2
\tag 2.6
$$
is nondecreasing  with decreasing $\Delta y_1$ as soon as $p_1\le p_2$ and 
$\Delta y_1/\Delta x_1\ge \Delta y_2/\Delta x_2$. This means 
that the maximum is  achieved at a convex $\psi$. 
Similarly, the expression (2.6) is nondecreasing with increasing 
$\Delta y_1$
if $p_1\ge \b$ and 
$\Delta y_1/\Delta x_1<\d/(1-\d)$, and therefore the maximum 
is achieved at a $\psi\in \Psi$. 

Next we note that 
$$
\F_k(\psi)\le \sum_{i=1}^k\z(\psi'(\a F(i/k)),i/k)\a \Delta F_k(i), 
$$
while 
$$
\F(\psi)\ge 
\sum_{i=1}^k \z(\psi'(\a F(i/k)),i/k)\a \Delta F_{k+1}(i). 
$$
Therefore, $\max_\psi \F_k(\psi)\le \max_\psi \F(\psi) +\O(1/k)$. 
As a lower bound is obtained similarly, this concludes the 
proof for nice distribution functions $F$.

To prove the general case, we again use Lemmas 2.1 and 2.4. 
For an arbitrary distribution function, choose nice $F_N^{\pm}$ 
so that $F_N^-\le F$ and $F\cdot 1_{(1/N, 1]}\le F_N^+$ and 
$F_N^{\pm}\to F$ and $b(F_N^{\pm})\to b(F)$. 
Then $c(\a,F_N^{\pm})\to c(\a,F)$. By Lemma 2.1, 
it immediately follows that
$\limsup H/m \le c(\a, F)$ a.s. 

The lower bound, however, does not immediately follow 
as $F$ is not below $F_N^+$. The remedy for this is 
to assume that $F(1/N)<1/2$, 
replace $\a$ with $\a'<\a$, and observe that the 
distribution $F$ will induce, with 
probability exponentially close to 1,
at least $(\a-\a')m/4$ probabilities $p_j\ge 1/N$.
Therefore the length of the longest increasing path 
in a $m\times \a m$ 
matrix using $F$ is eventually above the length of the 
longest increasing path 
in a $m\times \a' m$ 
matrix using $F_N^+$. By Lemma 2.4, 
$\liminf H/m \ge c(\a, F)$ a.s. $\square$\enddemo
 
\demo{Proof of Theorem 1 when $\a\le \a_c'$}
 
Applying the same strategy as before we
construct sequences
$\{F_N^{\pm}\}$ of distribution functions which satisfy 
$F_N^-\le F\le
F_N^+$
and for which Theorem 1 already holds, and such that
$c(F^-_N)$ and $c(F^+_N)$ approach the same 
limit as $N\to\infty$. Lemma 2.1 will then complete the 
proof. (We suppress $\a$ from the notation, since it 
is the same throughout this proof.) 

Take a sequence $\eN\searrow0$ such that $b-\eN$ are points of
continuity of $F$. Let $F^{\pm}_N$ agree with $F$ outside $[b-\eN, b)$, 
while 
on $[b-\eN, b)$ the two functions are constant: $F_N^-\equiv 
F(b-\eN)$ and  $F_N^+\equiv 1$. Let $\e_N=1-F(b-\eN)$; note 
that $\e_N\to 0$ and 
$dF_N^-=1_{(0,\,b-\eN)}\,dF+\e_N\,\delta_b$ and 
$dF_N^+=1_{(0,\,b-\eN)}\,dF+\e_N\,\delta_{b-\eN}$. 
Clearly  the already proved 
part of Theorem 1 applies to both $F_N^+$ and $F_N^-$. 

We proceed to show that $a(F_N^-)\to b$. If
this does not hold, the fact that   
$a(F_N^-)>b(F_N^-)=b$ implies that there exists 
an $\eta>0$ so that $a(F_N^-)\ge b+\eta$
along a subsequence. Then $\d=\<p(1-p)
\[(b-p)^{-2}-(b+\eta-p)^{-2}\]\>>0$ and 
$$
\align 
1&=\a\int_0^{b-\eN}\frac {x(1-x)}{(a(F_N^-)-x)^2}\,dF+
\a\,\e_N\,\frac {b(1-b)}{(a(F_N^-)-b)^2}\tag 2.7\\
&\le \a\int_0^b\frac {x(1-x)}{(b+\eta-x)^2}\,dF
+\a\,\e_N\,\frac {b(1-b)}{\eta^2}\\
&\le -\d+\frac {\a}{\a_c'}+\a\,\e_N\,\frac {b(1-b)}{\eta^2}, 
\endalign 
$$
along the same subsequence. 
As $N\to\infty$, this yields a contradiction with $\a\le \a_c'$. 
 
Now 
$$
c(F_N^-)=a(F_N^-)+\a\,(1-a(F_N^-))\(
\int_0^{b-\eN}\frac{x}{a(F_N^-)-x}\,dF+\e_N
\frac{b}{a(F_N^-)-b}\). 
\tag 2.8
$$
By (2.7), 
$$
\e_N\frac{b}{a(F_N^-)-b}\le \frac{a(F_N^-)-b}{\a (1-b)} \to 0. 
$$
To show that 
$$
\<1_{\{p\le b-\eN\}}p/(a(F_N^-)-p)\>\to \<p/(b-p)\>
\tag 2.9
$$
we note that the integrand on the left of (2.9) is uniformly 
integrable (as it is bounded by 
$p/(b-p)$, which is square--integrable) 
and converges to the 
integrand on the right a.s. 
By (2.8) and (2.9), 
$$
c(F_N^-)\to b+\a(1-b)\<p/(1-p)\>. 
$$
The argument for $c(F_N^+)$ is very similar and hence omitted. 
$\square$
\enddemo

\demo{Proof of Theorem 1 when $\a\ge \a_c$}
If $\a\nearrow\a_c$, then $a(\a,F)\nearrow 1$ and hence 
$c(\a,F)\nearrow 1$. 
$\square$
\enddemo

We note that the above proof of Theorem 1 actually 
shows exponential convergence to $c$, that is, 
(2.1) in Lemma 2.3 holds in random environment as well.
Also, once probabilities are ordered using Lemma 2.1, 
one could investigate convergence, in the sense of 
[DZ] and [Sep1], of a longest increasing 
path in $A$ to the maximizer of $\F(\psi)$. 
This is easy to prove if $\F$ is nice 
(cf. Lemma 2.5), but it actually holds whenever 
the maximizer is unique.  

We conclude this section by showing that the 
deterministic case indeed has no fluctuations. 

\proclaim{Proposition 2.6} Assume that $b<1$ and $\a> \a_c$. 
Then $P(H=m)$ converges to 1 exponentially fast
(and therefore $P(H=m$ eventually$)=1$). 
\endproclaim

\demo{Proof}  We begin by modifying the construction from
Section 3.3.1 of [GTW1]. Recall that random $m\times n$ matrix is the 
lower left corner of an infinite random matrix. 
For an $(i,j)\in \bN\times\bN$, let 
$\eta_{(i,j)}=\inf\{k\ge 1:  \e_{(i+k, j)}=0\}$ be the relative 
position of the first 0 above $(i,j)$ and 
$\xi_{(i,j)}=\inf\{k\ge 1: \e_{(i, j+k)}=1\}$
the relative position of the first 1 to the right of $(i,j)$.


Now define i.i.d. two--dimensional random vectors  
$X_1=(\xi_1,\eta_1)$, $X_2=(\xi_2, \eta_2),\dots$ as follows: 
$$
\alignat2 
&\xi_1=\xi_{(0,1)}, &&\quad \eta_1= \eta_{(\xi_1,1)},\\ 
&\xi_2=\xi_{(\xi_1,1+\eta_1)}, &&\quad\eta_2=\eta_{(\xi_1+\xi_2,1+\eta_1)},\\
&\xi_3=\xi_{(\xi_1+\xi_2,1+\eta_1+\eta_2)}, 
&&\quad\eta_2=\eta_{(\xi_1+\xi_2+\xi_3,1+\eta_1+\eta_2)},\\
&\dots
\endalignat 
$$
Let $S_k=(0,1)+X_1+\dots +X_k$ be the 
corresponding random walk, and $T_m$ (resp. $T_n'$) be the first 
time $S_k$ is in $\{(x,y): x>n\}$ (resp.  $\{(x,y): y>m\}$). 
If $T_m'<T_n$ then there is an increasing path of 1's 
inside the $m\times n$ rectangle which goes through 
its ``roof'' without skipping a row, thus 
$$
\{H<m\}\subset\{T_n\le T_m'\}.
$$
Therefore, we need to show that $P(T_n\le T_m')$ 
goes to 0 exponentially fast. To this end, note that, 
for any $\e>0$,
$$
P(T_n\le T_m')\le P(T_n\wedge T_m'\le \e m)+\sum_{k=\e m}^\infty 
P(\xi_1+\dots +\xi_k\ge \a (\eta_1+\dots +\eta_k)).
\tag 2.10
$$
If we show that $\xi_1$ and $\eta_1$ have exponential tails, 
and that $E(\xi_1)-\a E(\eta_1)<0$, the we can choose 
a small enough $\e>0$ so that the upper bound in (2.10) decays exponentially.
First, $P(\xi_1\ge k)=\<1-p\>^{k-1}$ and so $E(\xi_1)=1/\<p\>$. 
Moreover, the conditional distribution of $p$ given that 
a single coin flip gives 1 is 
$$
dF_1(x)=\frac 1{\<p\>}x\, dF(x), 
$$
therefore
$$
P(\eta_1\ge k)=\int_0^1 x^{k-1}\, dF_1(x)= \frac {\<p^k\>}{\<p\>}, 
$$
and so $E(\eta_1)=\<p/(1-p)\>/\<p\>$. $\square$
\enddemo

\subheading{3. The saddle point method and fluctuations}

Throughout this section, we assume that $b<1$ and that  
$\al=n/m$ is fixed
(but see Remark 3 at the end). In addition, our standing 
assumption will be that $$\a_c'<\a<\a_c.$$
We will investigate the limiting behavior of $P(H\le h)$ 
without using results proved in Section 2. An asymptotic analysis 
of this quantity when $\a<\a_c'$ is carried out in [GTW2]. 

We begin with deterministic inhomogeneous ODB, in which the $j$th column
is assigned a fixed {\it deterministic\/}
probability $p_j$. At first, our derivation 
will use a fixed $n$ and no particular properties of 
the eventual random choice of the environment. For notational 
convenience, we therefore drop the subscript $n$, which 
practically every quantity would otherwise have. 
See the discussion preceding the key formula (3.6), where the 
random environment is reintroduced. 

As explained in [GTW1], Sec. 2.2, we have
$$\Pr(H\le h)=\prod(1-p_j)^m\,D_h(\ph),$$
where $D_h$ is the $h\times h$ Toeplitz determinant with symbol
$$(1-z\inv)^{-m}\,\prod_{j=1}^n(1+r_jz)$$
and $r_j=p_j/(1-p_j)$. Applying an identity of Borodin and Okounkov
([BO], see also [BW])
this becomes
$$\Pr(H\le h)=\det\,(I-K_h),$$
where $K_h$ is the infinite matrix acting on $\ell^2({\bZ}^+)$ with
$j,k$ entry
$$K_h(j,k)=\sum_{\l=0}^{\iy}(\ph_-/\ph_+)_{h+j+\l+1}\;
(\ph_+/\ph_-)_{-h-k-\l-1}.$$
The subscripts here denote Fourier coefficients and the functions
$\ph_{\pm}$ are
the Wiener-Hopf factors of $\ph$, so
$$ \ph_+(z)=\prod_{j=1}^n(1+r_jz),\ \ \ \ph_-(z)=(1-z\inv)^{-m}.$$
The matrix $K_h$ is the product of two matrices, with $j,k$ entries give
by
$$(\ph_+/\ph_-)_{-h-j-k-1}={1\ov 2\pi i}
\int\prod(1+r_jz)\;(z-1)^m\,z^{-m+h+j+k}\,dz$$
and
$$(\ph_-/\ph_+)_{h+j+k+1}={1\ov 2\pi i}
\int\prod(1+r_jz)\inv\;(z-1)^{-m}\,z^{m-h-j-k-2}\,dz.$$
The contours for both integrals go around the origin once
counterclockwise; in the second
integral 1 is on the inside and all the $-r_j\inv$ are on the outside.

Eventually we let $m,\,n\ra\iy$ and will take $h=cm+sm^{1/3}$ where $c$,
as yet to be
determined, gives the transition between the limiting probability being
zero and
the limiting probability being one. In [GTW1] we considered the case
where
all the $p_j$ were the same. We found that with $c$ chosen as
in Lemma 2.3 we could do a
steepest descent analysis. The conclusion 
was that the product of the two matrices
scaled,
by means of the scaling $j\ra m^{1/3}x,\ k\ra m^{1/3}y$, to the square
of the integral
operator on $(0,\iy)$ with kernel ${\rmAi}(gs+x+y)$, where $g$ is
another
explicitly determined constant. This gave the limiting result
$$\lim_{n\ra\iy}\,\Pr(H\le cm+sm^{1/3})=F_2(g s),$$
where $F_2(s)$ is the Fredholm determinant of the Airy kernel on
$(s,\,\iy)$.
We can do very much the same here. If $h=cm+sm^{1/3}$ and we set
$$\ps(z)=\prod(1+r_jz)\;(z-1)^m\,z^{-(1-c)\,m}$$
then
$$(\ph_+/\ph_-)_{-h-j-k-1}={1\ov 2\pi i}\int
\ps(z)\,z^{s\,m^{1/3}+j+k}\,dz,\tag{3.1}$$
$$(\ph_-/\ph_+)_{h+j+k+1}={1\ov 2\pi
i}\int\ps(z)\inv\,z^{-s\,m^{1/3}-j-k-2}\,dz.
\tag{3.2}$$
To do an eventual steepest descent we define
$$\si(z)={1\ov m}
\log\,\ps(z)={\al\ov
n}\,\sum_{j=1}^n\log\,(1+r_jz)+\log\,(z-1)+(c-1)\log\,z,
$$
and look for zeros of
$$\si'(z)={\al\ov n}\,\sum_{j=1}^n{r_j\ov 1+r_jz}+{1\ov z-1}+{c-1\ov
z}.\tag{3.3}$$
The number of zeros equals one plus the number of distinct $r_j$. There
is a zero
between two consecutive $1/r_j$ and, in general, two other zeros which
which are either
unequal reals or a pair of complex conjugates. In the exceptional case
there
is a single real zero of multiplicity two. We choose $c$ so that we are
in this exceptional
case. If the double zero is at $z=u$ then $u$ and $c$ must
satisfy the pair of equations
$${\al\ov n}\,\sum_{j=1}^n{r_j\ov 1+r_ju}+{1\ov u-1}+{c-1\ov u}=0,$$
$${\al\ov n}\,\sum_{j=1}^n\left({r_j\ov 1+r_ju}\right)^2+{1\ov
(u-1)^2}+{c-1\ov u^2}=0.$$
If we multiply the second equation by $u$ and subtract we get
$${\al\ov n}\,\sum_{j=1}^n{r_j\ov
(1+r_ju)^2}={1\ov(u-1)^2}.\tag{3.4}$$
The first equation gives
$$c={1\ov 1-u}-{\al\ov n}\sum_{j=1}^n{r_ju\ov 1+r_ju}.\tag{3.5}$$
Conversely, if the second pair of equations is satisfied then so is the
first. 

\proclaim {Lemma 3.1} Assume that $\al n\inv\sum r_j<1$ and set $\ub=\max
\{-1/r_j\}$.
Then equation ({3.4}) has a unique solution $u\in (\ub,\,0)$ and if
$c$ is then
defined by ({3.5}) we have $c\in (0,\,1)$.
\endproclaim

\demo{Proof} The left side of ({3.4}) decreases from $\iy$ to
$\al n\inv\sum r_j$
as $u$ runs over the interval $(\ub,\,0)$ whereas the right side
increases and has the
value 1 at $u=0$. Our assumption implies the first statement of the
lemma. As for the
second, $c>0$ since $u<0$ and each $1+r_ju>0$. Moreover, 
Schwarz's inequality,
our assumption and ({3.4}) give
$${\al\ov n}\sum{r_j\ov 1+r_ju}\le \left\{{\al\ov n}\sum
r_j\right\}^{1/2}\,
\left\{{\al\ov n}\sum {r_j\ov (1+r_ju)^2}\right\}^{1/2}<{1\ov1-u}.$$
Hence
$$c<{1\ov 1-u}-{u\ov 1-u}=1.$$
$\square$
\enddemo

To derive the asymptotics using steepest descent we have to compute
$\si'''(u)$ and
understand the steepest descent curves. For the first we multiply
({3.3}) by $z$,
differentiate twice and use the fact that $\si'(u)=\si''(u)=0$ to obtain
$$u\,\si'''(u)=-{2\al\ov n}\,\sum_{j=1}^n{r_j^2\ov (1+r_ju)^3}+{2\ov
(u-1)^3}.$$
Note that $\si'''(u)>0$ since $u<0$.

There are three curves emanating from $z=u$ on each of which $\Im\,\si$
is constant. One is
$\Im\,z=0$, which is of no interest. The other two come into $u$ at
angles $\pm\pi/3$
and $\pm2\pi/3$. Call the former $C^+$ and the latter $C^-$.
Approximate shapes of these curves are illustrated in 
Figure 4.  
For the integral involving $\ps(z)$ we want $|\ps(z)|$ to have a maximum
at the point $u$ on 
the curve and for
the integral involving $\ps(z)\inv$ we want $|\ps(z)|$ to have a minimum
at $u$. Since $\si'''(u)>0$ the curve for $\ps(z)$ must be $C^+$ and the
curve for
$\ps(z)\inv$ must be $C^-$.
 
As for the global natures of the curves, $C^{\pm}$ can only end at a
zero of $\ps(z)^{\pm1}$,
at a zero of $\si'(z)$, or at 
infinity. The two curves are simple and
cannot intersect
since $|\ps(z)|$ is
decreasing on $C^+$ as we move away from $z=u$ while $|\ps(z)|$ is
increasing on $C^-$. It follows that $C^+$ closes at $z=1$, while the two
branches of
$C^-$ go to infinity. From the fact that
$$\Im\,\si(z)={\al\ov
n}\,\sum_{j=1}^n\arg\,(1+r_jz)+\arg\,(z-1)+(c-1)\arg\,z$$
is constant on $C^-$ we can see that the two branches go to infinity in
the directions
$\arg z=\pm c\pi/(c+\al(1-\nu))$, where $\nu$ is the fraction of $r_j$
equal to zero
(which is the same as 
fraction of the $p_j$ equal to zero). Observe that in the integral in
({3.1}) the
path can be deformed into $C^+$ and in the integral in ({3.2}) the
path
can be deformed into $C^-$. Both contours will be described downward
near $u$.

\smallskip
\epsfysize=0.3\vsize
\centerline{\epsffile{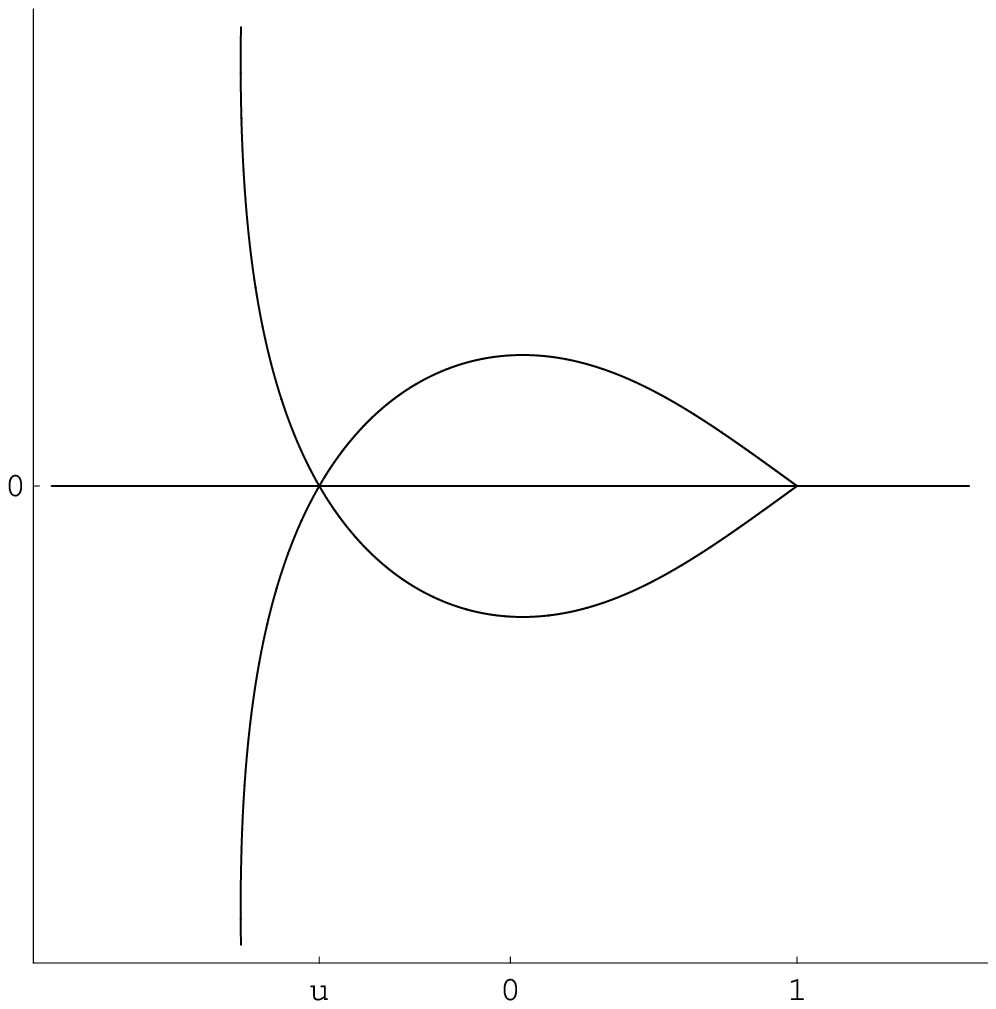}}
\smallskip
 \centerline {Figure 4. The steepest descent curves $C^{\pm}$ as
described in the text.}
\smallskip

To see formally what steepest descent gives, 
we replace our matrices
$M(j,k)$ depending
on the parameter $m$ and acting on $\ell^2({\bZ}^+)$ by kernels
$m^{1/3}M(m^{1/3}x,\,m^{1/3}y)$ acting on $L^2(0,\,\iy)$. Thus
({3.1}) becomes the
operator with kernel
$${1\ov 2\pi i}m^{1/3}\int e^{m\si(z)}\,z^{m^{1/3}(s+x+y)}\,dz.$$

If steepest descent worked, the main contribution would come from the
immediate neighborhood
of $z=u$. We would set $z=u+\z$, make the replacements
$$\si(z)\ra\si(u)+{1\ov6}\si'''(u)\z^3,\ \ \ z\ra u\,e^{\z/u},$$
in the integral and integrate (downwards) on the rays $\arg\z=\pm\pi/3$.
The above integral becomes
$$e^{m\si(u)}\,u^{m^{1/3}(s+x+y)}{1\ov 2\pi i}m^{1/3}\int
e^{{m\ov6}\si'''(u)\z^3+m^{1/3}(s+x+y)\z/u}\,d\z,$$
and we can then replace the rays by the imaginary axis (downwards). The
variable change
$\z\ra-i\z/m^{1/3}$ replaces this by
\footnote{Recall that the Airy function is defined by
${\rmAi}(x)={1\over 2\pi} 
\int_{-\iy}^{\iy} e^{i\z^3/3+i x \z}\, d\zeta.$}
$$-e^{m\si(u)}\,u^{m^{1/3}x}{1\ov 2\pi}\int_{-\iy}^{\iy}
e^{{i\ov6}\si'''(u)\z^3-i(s+x+y)\z/u}\,d\z$$
$$=-e^{m\si(u)}\,u^{m^{1/3}(s+x+y)}\,|u|g\,{\rmAi}(g(s+x+y)),
$$
where we have set
$$g=|u|\inv\left\{{1\ov2}\si'''(u)\right\}^{-1/3}.$$
Thus, if we multiply the matrix entries on the left
side of ({3.1})
by
$$-e^{-m\si(u)}\,u^{-m^{1/3}s-j-k},$$
then the result has the scaling limit the operator
on $L^2(0,\,\iy)$ with kernel
$$|u|g\,{\rmAi}(g(s+x+y)).$$
Similarly if we multiply the
matrix entries on the left side of ({3.2}) by
$$-e^{m\si(u)}\,u^{m^{1/3}s+j+k},$$
then the result has the scaling limit the operator
on $L^2(0,\,\iy)$ with kernel
$$|u|\inv g{\rmAi}(g(s+x+y)).$$
It follows that
the product of the two matrices has in the limit the same Fredholm
determinant as
the operator with kernel
$$g^2\int_0^{\iy}{\rmAi}(g(s+x+z))\,{\rmAi}(g(s+z+y))\,dz$$
$$=g\int_0^{\iy}{\rmAi}(g(s+x)+z)\,{\rmAi}(g(s+y))\,dz$$
which in turn has the same Fredholm determinant as the kernel
$$\int_0^{\iy}{\rmAi}(gs+x+z)\,{\rmAi}(gs+z+y)\,dz.$$
This Fredholm determinant equals $F_2(gs)$.

Assuming the argument we sketched above 
goes through we will have shown that, in some sense,
$$\lim_{n\ra\iy}\,\Pr(H\le cm+sm^{1/3})=F_2(g s),$$
where $c$ and $g$ are as above and determined once we know the $p_j$ and
$\al$.

We begin the rigorous justification by introducing some 
notation. Recall that we consider a random environment
in which the probabilities $p_j$ are
chosen
independently with distribution function $F$. 
We explained the notation $H_n$ after Lemma 2.2; in addition, we 
give the subscript $n$ to the 
quantities $\sig_n(z)$, $u_n$, $g_n$, and curves  
$C_n^{\pm}$ to emphasize that they are 
functions of $p_1,\dots, p_n$. 
Therefore 
$$\Pr(H\le h)=\lan\Pr(H_n\le h)\ran,$$
where $\lan\,\cdot\,\ran$ is the expected value
with respect to $p_1,\dots, p_n$. 

Our object is to show that
with
probability one, for each fixed $s$, 
$$\Pr(H_n\le c_n\,m+s\,m^{1/3})=F_2(g_n \,s)+o(1)
\tag{3.6}$$
as $n\ra\iy$. 
We will demonstrate these asymptotics by pointing out the necessary 
modifications to the argument in [GTW1]. 

All the $-1/r_j$ in our previous discussion are
contained in the interval $(-\iy,\,\xi]$, where $\xi=1-1/b$. (Recall
that $b$
is the maximum of the support of $dF$.) Let $F_n$ be the empirical
distribution function given by 
$$dF_n=n\inv\sum\dl_{p_j}$$
and let $\<\,\cdot\,\>_{F_n}$
denote the integration with respect to $dF_n$.  
Recall the
Glivenko-Cantelli
theorem, which says that, with probability one, 
$F_n$ converges uniformly to $F$
as
$n\ra\iy$. 

We first show that, under our standing assumptions, the quantities
$c_n$ and $u_n$ of Lemma 3.1
converge almost surely as $n\ra\iy$ to the corresponding quantities 
associated with the distribution function $F$. Recall that we set
$r=p/(1-p),\; p=r/(1+r)$. We remark that $c_0$ in the following 
lemma is the same as $c$ in Theorem 1, and 
$u_0= (a-1)/a$. The notation has changed 
to conform with (3.4) and (3.5), which are 
in turn chosen to connect with the saddle--point approach in [GTW1]. 

\proclaim{Lemma 3.2}  The equation
$$\al\lan{r\ov (1+ru_0)^2}\ran={1\ov(u_0-1)^2}$$
has a unique solution $u_0\in (\xi,\,0)$ and if $c_0$ is then defined by
$$c_0={1\ov 1-u_0}-\al\lan{ru_0\ov 1+ru_0}\ran$$
we have $c_0\in (0,\,1)$.
\endproclaim

\demo{Proof} The argument goes almost exactly as for Lemma 3.1. The 
assumption  $\a<\a_c$ is equivalent to
$\al\lan r/(1+ru)^2\ran>1/(u-1)^2$ when $u=\xi$, while 
$\a>\a_c'$ yields the opposite inequality when $u=0$. $\square$\enddemo

Note that one obtains 
$u_n$ as $u_0$, except that the expectation $\<\,\cdot\,\>$ is replaced
by the expectation $\<\,\cdot\,\>_{F_n}$.
 
\proclaim{Lemma 3.3} Almost surely, 
$u_n\ra u_0$ and $c_n\ra c_0$ as $n\ra\iy$.
\endproclaim

\demo{Proof} Integration by parts gives
$$\lan{r\ov (1+rz)^2}\ran=\int_0^{b/(1-b)}(1-F(p))\,{d\ov dr}{r\ov
(1+rz)^2}\,dr.$$
The derivative in the integrand is uniformly bounded for $z$ in any
compact subset of the
complement of $(-\iy,\,\xi]$. Hence the expected value is continuous in
$F$
and differentiable for $z\not\in(-\iy,\,\xi]$. Moreover,
$${\partial\ov\partial z}\left(\al\lan{r\ov
(1+rz)^2}\ran-{1\ov(z-1)^2}\right)$$
is negative, hence nonzero, at $z=u_0$. The statement concerning $u_n$
therefore
follows from the fact that $F_n\ra F$ uniformly and the implicit
function theorem. The
assertion for $c_n$ then follows by a similar 
integration by parts. $\square$\enddemo

\proclaim{Lemma 3.4}  There exists a
(deterministic) wedge $W$ with vertex $v>\xi$,
bisected by
the real axis to the left of $v$, such that, 
almost surely, the curves $C_n^{\pm}$
lie outside $W$ for sufficiently large $n$.\endproclaim

\demo{Proof} First we show that, if $\ve$ is small enough, $C_n^-$
is disjoint from the disc
$$D(\xi,\,\ve)=\{z\,:\,|z-\xi|\le\ve\}.$$
From the facts that
$\si_n'(u_n)=\si_n''(u_n)=0,\ \si_n'''(u_n)>0$, and $\si_n'(z)\ne 0$
for $z\in(\xi,\,u_n)$,
it follows that $\si_n$ is strictly increasing in the interval $(\xi,\, u_n)$. 
Therefore, we can choose small enough $\e>0$ and 
$\dl>0$ so that $\si_n(\xi+\ve)<\si_n(u_n)-2\dl$ 
for all large enough $n$. 
In
addition, if $\ve$
is small enough,
$$\log|z-1|+(c_n-1)\log|z|<\log|\xi+\ve-1|+(c_n-1)\log|\xi+\ve|+\dl$$
for all $z\in D(\xi,\,\ve)$.
Now each $|1+r_jz|$, and so its logarithm, achieves its maximum on
$D(\xi,\,\ve)$ at the point
$z=\xi+\ve$. 
By combining the last three observations, 
we see that everywhere on $D(\xi,\,\ve)$ we have
$$\Re\,\si_n(z)<\si_n(\xi+\ve)+\dl<\si_n(u_n)-\dl.$$
Since $\Re\,\si$ achieves its minimum on $C_n^-$ at $z=u_n$, the curve
must be
disjoint from $D(\xi,\,\ve)$.


\epsfysize=0.9\vsize
\centerline{\epsffile{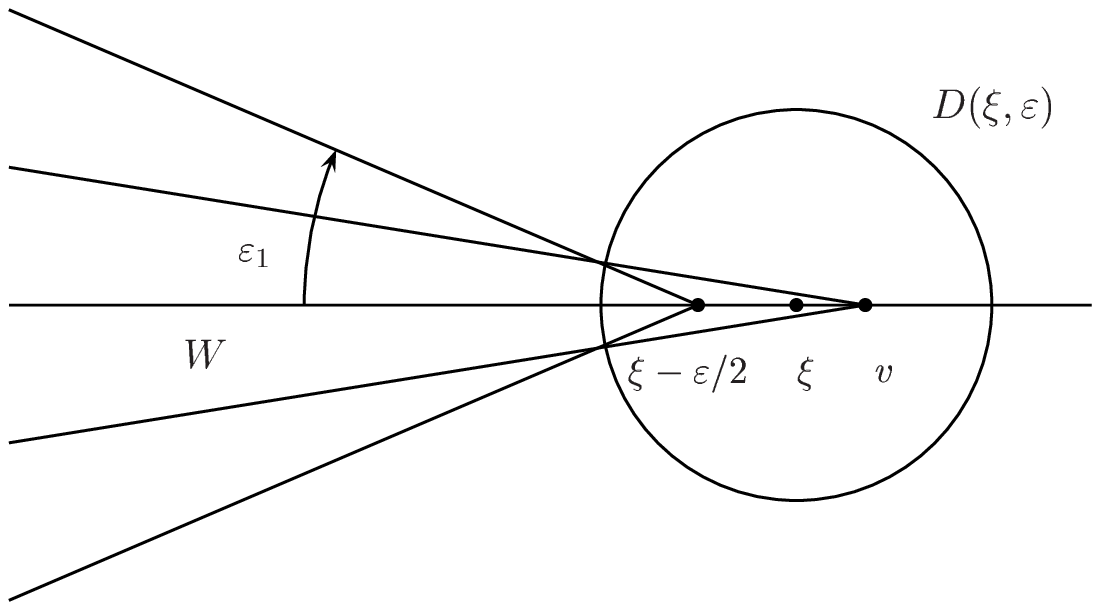}}
\vskip-10cm
\centerline {Figure 5. Wedge $W$,  angle $\ve_1$, and disk $D(\xi,\ve)$ 
as described in the proof of Lemma 3.4.}
\smallskip

\vskip0cm

For a small $\e_1>0$ (possibly much smaller than $\e$), 
denote by $W'$ the wedge with vertex $\xi-\ve/2$ bounded by the real
axis to the left of $\xi-\ve/2$ and the ray
$\arg\,(z-\xi+\ve/2)=\pi-\ve_1$.
Our next step is to 
show that $C_n^-$ is
disjoint
from  $W'$ if $\ve_1$ is small enough. 
As $\Im\, \sig_n(z)$ is constant 
on the portion of $C_n^-$ in the upper half-plane, 
$${\al\ov
n}\,\sum_{j=1}^n\arg\,(1+r_jz)+\arg\,(z-1)+(c_n-1)\arg\,z=c_n\pi,$$
where all arguments lie in $[0,\,\pi]$.  For $z\in W'$, 
$$\arg\,(z-1)+(c_n-1)\arg z\ge c_n\arg z\ge c_n(\pi-\ve_1).$$
Since $b$ is in the support of $dF$, 
the strong law implies that, almost surely,   
at least a positive fraction $\eta$ of the $-1/r_j$ lie in the interval
$[\xi-\ve/2,\,\xi]$ for large enough $n$. The contribution 
of these terms (and nonnegativity of the others) 
in the following sum 
provides a lower bound valid for $z\in W'$: 
$${\al\ov n}\,\sum_{j=1}^n\arg\,(1+r_jz)\ge \a \eta(\pi-\e_1).$$
Hence, for $z\in W'$, 
$$\Im\,\si_n(z)\ge\al\eta(\pi-\ve_1)+c_n(\pi-\ve_1)=
(\al\eta+c_n)\,(\pi-\ve_1)>c_n\pi,$$
if $\ve_1$ is chosen to be small enough. Therefore, $C_n^-$
is disjoint from
the wedge $W'$ for small enough $\ve_1$. By symmetry, $C_n^-$ is disjoint
from the
reflection of $W'$ over the imaginary axis. We have shown that the
curve is also
disjoint from $D(\xi,\,\ve)$ and the union of the disc and the two
wedges contains a
wedge of the form described in the statement of the lemma.

This establishes the statement of the lemma concerning $C_n^-$. Since
$C_n^+$ is to the
``right'' of $C_n^-$ (it begins to the right and they cannot cross), the
statement for
$C_n^+$ follows automatically. $\square$\enddemo

In the following lemma $\si_0$ denotes the function $\si$ associated
with the
distribution $F$,
$$\si_0(z)=\al\lan\log\,(1+rz)\ran+\log\,(z-1)+(c_0-1)\log\,z,$$
and 
$$
g_0=|u_0|^{-1}\left\{\frac 12 \sig_0'''(u_0)\right\}^{-1/3}. 
$$

\proclaim{Lemma 3.5} Almost surely, $z\,\si_n'(z)\ra z\,\si_0'(z)$ uniformly
outside the wedge $W$ of Lemma~3.4.\endproclaim

\demo{Proof} We have
$$z\si_n'(z)-z\si_0'(z)=\al\int{rz\ov 1+rz}d(F_n(p)-F(p))+ c_n-c_0$$
$$=\al\int_0^{b/(1-b)}(F_n(p)-F(p))\,{z\ov (1+rz)^2}\,dr+c_n-c_0.$$
The last term goes to 0 by Lemma 3.3.
The last factor in the integrand is uniformly bounded for
$r\in(0,\,b/(1-b))$,
$z\not\in W$ and $z$ bounded. Thus $z\,\si_n'(z)\ra z\,\si_0'(z)$
uniformly on
bounded subsets of the complement of $W$. If $z$ is sufficiently
large and outside $W$ then it is outside some wedge with vertex 0
bisected by the
negative real axis, and on the complement of any such wedge
$$\int_0^{b/(1-b)}{|z|\ov |1+rz|^2}\,dr$$
is uniformly bounded. Thus $z\,\si_n'(z)\ra z\,\si_0'(z)$
uniformly throughout the complement of $W$. $\square$\enddemo

The preceding lemmas show that the curves $C_n^{\pm}$ are uniformly
smooth, as we now
argue. The function $\si_0'(z)$ can have no other zero in the
complement of
$(-\iy,\,\xi)$
than at $z=u_0$. This follows from uniform convergence and the fact that
the
corresponding statement holds for the $\si_n(z)$. Thus the 
functions $\si_n'(z)$ are
uniformly bounded away from zero on compact subsets not containing
$u_0$. As we move outward
(i.e., away from $u_n$) along $C_n^{\pm}$,
$\Im\,\si$ is constant and $\Re\,\si$ is increasing on $C_n^-$ and
decreasing
on $C_n^+$. It follows that if $s$ measures arc length on the 
curves then, for $z\in C_n^{\pm}$, 
$${dz\ov ds}=\mp{|\si_n'(z)|\ov \si_n'(z)}.\tag{3.7}$$
This shows that the $C_n^{\pm}$ are uniformly smooth on compact sets (to
be more precise,
the portions in the upper and lower half-planes are). Moreover, 
they are
uniformly close
on compact sets to the corresponding curves $C_0^{\pm}$ for the
distribution
function $F$. In particular, the length of $C_n^+$ is $O(1)$. 

To see what happens for large $z$ on $C_n^-$ observe that
$$\lim_{z\ra\iy}z\,\si_n'(z)=\al+c_n>0$$
uniformly in $n$. This and ({3.7}) show that $|z|$ is increasing as
we
move far enough out along $C_n^-$. If $\Gamma$
is an arc of $C_n^-$ going from $a$ to $b$ then
$$\int_{\Gamma}|\si_n'(z)|\,ds=\int_{\Gamma}\si_n'(z)\,dz=\si_n(b)-\si_n(a).$$
Hence the length of $\Gamma$ is at most $|b-a|$ times
$${\max_{z\in[a,b]}|\si_n'(z)|\ov \min_{z\in\Gamma}|\si_n'(z)|},$$
where $[a,\,b]$ is the line segment joining $a$ and $b$, as long as this
segment
does not meet $(-\iy,\,\xi)$. It follows from the above, for example,
that
the $L^1$ norm of the function $(1+|z|^2)\inv$ on $C_n^-$ is $O(1)$.

In [GTW1] we needed asymptotics with
error bounds for all $j,k\le h$ and this required a more careful
analysis of the
integrals in ({3.1}) and ({3.2}) than we indicated; instead
of the steepest descent curves passing through the same point they pass
through
different, but nearby, points. With what we now know we can
show that these curves are uniformly smooth with uniformly regular
behavior near
infinity, and this is what is needed to see that in our case the
asymptotics hold uniformly
in $n$.

\proclaim{Lemma 3.6} Almost surely, $g_n\ra g_0\ne 0$ and ({3.6}) holds.
\endproclaim

\demo{Proof} The first statement follows from Lemmas 3.3 and 3.5 and
the fact that the
$\si_n'(z)$ have only two zeros outside $(-\iy,\,\xi]$ counting
multiplicity, and therefore
$\si'''(u_0)\ne0$. 

To establish (3.6), one now has to go through the  
the steepest descent argument in Sec. 3.1.2 of [GTW1],
and make some obvious changes, justified 
by the results of this section. For the analogue of Lemma 3.1 there, for
example, we would add 
the phrase ``and all sufficiently large $n$'' to the end of the
statement. At the
end of second sentence of the proof we would add the phrase ``since
$\si_n'''(u_n)$ is 
uniformly bounded away from zero the length of $C_n^+$ is $O(1)$.''
After the last sentence 
we would add
``again since the length of $C_n^+$ is $O(1)$.''
Analogous changes need to be made throughout the argument and 
we skip further details.
$\square$\enddemo

\demo{Proof of Theorem 3} By Lemma 3.6, we can take $G_n=c_nm$. 
$\square$\enddemo

\demo{Proof of Theorem 2}
Note first that 
$$\tau^2={\text Var}\left({ru_0\ov 1+ru_0}\right),$$
where $u_0$ is as in Lemma 3.2 (and, as we remarked earlier, $c_0=c$).  

The proof rests on the crucial property 
(3.6) and the fact that $n^{1/2}(F_n-F)$
converges in
distribution to a Brownian bridge $B$ with an appropriate 
covariance structure; in particular $B$ is a Gaussian random 
element in $D[0,1]$ ([Bil], Th. 14.3). By the Skorohod representation 
theorem, 
we can couple $F_n$ and $B$ on some 
probability space $\Omega_0$ so that 
$$n^{1/2}(F_n-F)\to B\tag 3.8$$
in fact converges for every $\omega\in \Omega_0$
([Bil], Theorem 6.7). We now prove that, under this 
coupling, the
solution $u_n$ of 
({3.4}) satisfies 
$$u_n= u_0+n^{-1/2}\,U+o(n^{-1/2}),\tag 3.9$$ 
for every $\omega$ and 
for some Gaussian random variable $U$.

To establish (3.9), define 
$$
\align 
&\th_n(u)=\al\lan{r\ov(1+ru)^2}\ran_{F_n}-{1\ov(u-1)^2},\\
&\th_0(u)=\al\lan{r\ov(1+ru)^2}\ran-{1\ov(u-1)^2}.
\endalign 
$$
By Lemma 3.6 and its proof, 
there exists a 
(deterministic) neighborhood $\U\subset \bC$ 
of $u_0$ 
in which, with probability 1, $u_n$ (resp. $u_0$) is for large $n$ 
the unique solution to $\th_n(u)=0$ 
(resp. $\th_0(u)=0$). 
Therefore we can choose a fixed contour $C$ in $\U$ such that 
$u_n$ and $u_0$ are given by
$$u_n={1\ov 2\pi i}\int_C{\th_n'(u)\ov\th_n(u)}\,u\,du, \quad
u_0={1\ov 2\pi i}\int_C{\th_0'(u)\ov\th_0(u)}\,u\,du.$$

By (3.8), we have, uniformly for $u\in C$, 
$$\th_n(u)=\th_0(u)+n^{-1/2}\al\lan{r\ov(1+ru)^2}\ran_B+o(n^{-1/2}).$$
Here $\lan\,\cdot\,\ran_B$ is the expectation with respect to 
$dB$, but by integration by parts (as in the proof of Lemma 3.3) 
we can 
make $B$ appear 
in the integrand. 
Therefore
$${\th_n'(u)\ov\th_n(u)}={\th_0'(u)\ov\th_0(u)}+
n^{-1/2}\al\,{d\ov du}{\lan r/(1+ru)^2\ran_B\ov\th_0(u)}+o(n^{-1/2}).$$
If we multiply this identity 
by $u/2\pi i$ and integrate over $C$ the left side becomes $u_n$,
the first term on 
the right becomes $u_0$ while the second 
term becomes $n^{-1/2}\,U$ where 
$$U=-\al\,{\lan r/(1+ru_0)^2\ran_B\ov\th_0'(u_0)}$$ 
is a Gaussian random variable. This proves (3.9). 

Let 
$$\ph_n(u)=
{1\ov 1-u}-\a\<\frac{r}{1+ru}\>_{F_n},$$
so that $c_n=\ph_n(u_n)$. We claim that 
$$c_n=\ph_n(u_0)+O(n^{-1}).\tag 3.10
$$
To see this, we use the fact that $\ph_n'(u_n)=0$ to write
$$c_n=\ph_n(u_n)=\ph_n(u_0)+(u_n-u_0)^2\int_0^1
t\ph_n''(t\,u_n+(1-t)\,u_0)\,dt.$$
Thus, (3.10) follows from (3.9) and the uniform 
boundedness of the $\phi_n''(u)$ near $u=u_0$.  

Now, by the central limit theorem, 
$$
\sqrt n\({1\ov n}\sum_{j=1}^n{r_ju_0\ov 1+r_ju_0}-\<\frac{ru_0}{1+ru_0}\>\)
$$ 
converges in distribution to a Gaussian random variable 
$X$ with mean 0 and variance $\tau^2$. Therefore,
$$
\sqrt n(c_n-c_0)\convd \a X.
\tag 3.11
$$
Finally, (3.6) implies that,  
for any $\dl>0$, 
$$
P(-\d m^{1/2}\le H-c_n m\le \d m^{1/2})\to 1.
\tag 3.12
$$
(In fact, (3.6) implies that the above statement holds
with probability 1 before the expectation with respect 
to $p_1,\dots, p_n$ is taken, that is, if $H$ is replaced by $H_n$.) 
It follows from (3.11) and (3.12)  
that 
$$
\(H-c_0\,m\)/\sqrt m\convd\sqrt\al\,X,
$$
which concludes the proof. 
$\square$\enddemo

\noi{\bf Remark 1}. We did not need the full force of (3.6)
for the above proof to go through. Instead, a much weaker
property (3.12) suffices. \sp

\noi{\bf Remark 2}. As mentioned in the Introduction, 
independence of $p_n$ is not necessary for the results of 
this section to hold. Indeed, one only needs Glivenko--Cantelli
theorem for convergence in probability of $H/m$ to the time constant, 
hence ergodicity of $p_1,p_2,\dots$ is enough. Furthermore, 
a strong enough mixing property of this sequence  is sufficient 
for a normal fluctuation result. This follows from 
Billingsley's results in Section 22 of the first (1968) edition of [Bil]. \sp 

\noi{\bf Remark 3}. We assumed that $\a=n/m$ is fixed, but 
the proof of Theorem 2 remains valid with $n=\a m+o(\sqrt m)$. 

\comment 

\noi{\bf Remark 2}. If we recall that $r=p/(1-p)$ and define $a$ by
$u_0=(a-1)/a$ then
$a$ satisfies
$$\al\lan{p(1-p)\ov (a-p)^2}\ran=1,$$
our assumptions become $\al'_c<\al<\al_c$ and we can see that $c_0$ is
equal to
$c(\al,F)$ as given in the the middle formula in the 
statement of Theorem 1
of the introduction.   

\endcomment

\vskip1cm

\centerline{REFERENCES}

\vskip0.5cm

\item{\bf [Ale]} K.~S.~Alexander,  {\it
Approximation of subadditive functions and convergence rates in limiting-shape
results.\/} Ann.\ Probab.\ 25 (1997), 30--55.

\item{\bf [BDJ]}
J.~Baik, P.~Deift, K.~Johansson, {\it 
On the distribution of the length of the 
longest increasing subsequence
of random permutations.\/} J.\ Amer.\ Math.\ Soc. 12 (1999),
1119--1178.

\item{\bf [Bil]} P.~Billingsley, ``Convergence of Probability Measures.'' 
John Wiley, 1999.

\item{\bf [BO]} A.~Borodin, A.~Okounkov, {\it A Fredholm
determinant formula for Toeplitz determinants.\/}  
Int.\ Eqns.\ Operator Theo.\ 37 (2000), 386--396. 

\item{\bf [BR]} 
J.~Baik, E.~M.~Rains, {\it Limiting distributions for a polynuclear
growth model with external sources.\/}
J.\ Statist.\ Phys. 100 (2000), 523--541.

\item{\bf [BW]} E.~L.~Basor, H.~Widom, {\it 
On a Toeplitz determinant identity
of Borodin and Okounkov.\/}
Int.\ Eqns.\ Operator Theo.\ 37 (2000), 397--401.

\item{\bf [Dur]} R.~Durrett, ``Lecture Notes on 
Particle Systems and Percolation.'' Brooks/Cole, 1988.

\item{\bf [DZ]} J.--D.~Deuschel, O.~Zeitouni, {\it 
Limiting curves for i.i.d.\ records.\/} Ann.\ Probab.\ 23
(1995), 852--878.

\item{\bf [Gri]} D.~Griffeath, {\it Primordial 
Soup Kitchen.\/} {\tt psoup.math.wisc.edu}

\item{\bf [Gra]} J.~Gravner, {\it Recurrent ring 
dynamics in two--dimensional excitable cellular 
automata.\/} J.\ Appl.\ Prob.\ 36 (1999), 492--511. 

\item{\bf [GTW1]} J.~Gravner, C.~A.~Tracy, H.~Widom, 
{\it Limit theorems for height fluctuations in a class of discrete space
and time growth models. \/} 
J.\ Statist.\ Phys.\ 102 (2001), 1085--1132.

\item{\bf [GTW2]} J.~Gravner, C.~A.~Tracy, H.~Widom, 
{\it  A growth model in a random environment II: The composite regime.\/} 
In preparation.

\item{\bf [HW]} J.~M.~Hammersley, D.~J.~Welsh,  {\it 
First-passage percolation, subadditive processes, stochastic networks,
and generalized renewal theory.\/} In ``Bernoulli, Bayes, 
Laplace Anniversary Volume,'' J.~Neyman and L.~LeCam, editors, 
Springer-Verlag, 1965. Pages 61--110.

\item{\bf [ITW1]} A.~R.~Its, C.~A.~Tracy, H.~Widom,
{\it Random Words, Toeplitz Determinants and Integrable Systems. I.\/}
In ``Random Matrix Models and their Applications,'' Math.\ Sci.\
Res.\ Inst.\ Publications, Vol. 40,  P.~Bleher and A.~R.~Its, editors, 
Cambridge University Press, New York, 2001. Pages 245--258.

\item{\bf [ITW2]} A.~R.~Its, C.~A.~Tracy, H.~Widom,
{\it Random Words, Toeplitz Determinants and Integrable Systems. II.\/}
Physica D 152--153 (2001), 1085--1132. 

\item{\bf [Joh1]}
K.~Johansson, {\it Shape fluctuations and random matrices.\/} 
Commun.\ Math.\ Phys.\ 209 (2000), 437--476.

\item{\bf [Joh2]} K.~Johansson, {\it 
Discrete orthogonal polynomial ensembles and the Plancherel
measure.\/} Ann.\ Math.\ 153 (2001) 259--296. 

\item{\bf [Mea]} P.~Meakin, ``Fractals, scaling and growth 
far from equilibrium.'' 
Cambridge University Press, Cambridge, 1998.

\item{\bf [NS]} C.~M.~Newman, D.~L.~Stein,   
{\it Equilibrium pure states and nonequilibrium chaos.\/} 
J.\ Statist.\ Phys.\ 94 (1999), 709--722.

\item{\bf [NV]} C.~M.~Newman, S.~B.~Volchan, {\it 
Persistent survival of one-dimensional contact 
processes in random environments. \/} 
Ann.\ Probab.\ 24 (1996), 411--421. 

\item{\bf [PS1]} M.~Pr\"ahofer, H.~Spohn, {\it Universal distribution 
for growth processes in $1+1$ dimensions and random 
matrices\/}. Phys.\ Rev.\ Lett.\ 84 (2000), 
4882--4885.

\item{\bf [PS2]} 
M.~Pr\"ahofer, H.~Spohn, {\it Scale Invariance of the PNG Droplet
and the Airy Process.\/}  Preprint (ArXiv: math.PR/0105240).

\item{\bf [Rai]} E.~M.~Rains, {\it A mean identity 
for longest increasing subsequence problems.\/}
Pre\-print (arXiv: math.CO/0004082).

\item{\bf [Sep1]} T.~Sepp\" al\" ainen, {\it Increasing sequences of 
independent points on the planar lattice.\/} 
Ann.\ Appl.\ Probab. 7 (1997), 886--898.
 
\item{\bf [Sep2]} T.~Sepp\"al\"ainen, {\it 
Exact limiting shape for a simplified model of first-passage
percolation on the plane.\/} Ann.\ Probab.\ 26 (1998), 1232--1250.

\item{\bf [SK]} T.~Sepp\"al\"ainen,  J.~Krug, {\it Hydrodynamics and 
platoon formation for a  totally asymmetric
exclusion model with particlewise disorder.\/} 
J.\ Statist.\ Phys.\ 95 (1999), 525--567.

\item{\bf [Tal]} M.~Talagrand, 
{\it Huge random structures and mean field models for spin glasses.\/}
Doc.\ Math., Extra Vol. I (1998), 507--536. 

\item{\bf [TW1]} C.~A.~Tracy, H.~Widom, 
{\it Level spacing distributions and the Airy kernel.\/}
Commun.\ Math.\ Phys.\ 159 (1994), 151--174.

\item{\bf [TW2]} C.~A.~Tracy, H.~Widom,
{\it Universality of the Distribution Functions of 
Random Matrix Theory. II.\/} In 
``Integrable Systems: From Classical to Quantum,'' 
J.~Harnad, G.~Sabidussi and P.~Winternitz, editors, 
American Mathematical Society, Providence, 2000. Pages 251--264. 

\item{\bf [TW3]} C.~A.~Tracy, H.~Widom, {\it 
On the distributions of the lengths of the longest monotone 
subsequences in random words.\/} 
Prob.\ Theory Rel.\ Fields 
119 (2001), 350--380. 

\enddocument